\documentclass[12pt]{article}
\usepackage[margin=2cm]{geometry}
\usepackage{amsfonts}
\usepackage{amssymb}
\usepackage{amsthm}
\usepackage{amsmath}
\usepackage{latexsym}
\usepackage{color}

\begin{document}
\newtheorem{lemma}{Lemma}
\newtheorem{pron}{Proposition}
\newtheorem{thm}{Theorem}
\newtheorem{Corol}{Corollary}
\newtheorem{exam}{Example}
\newtheorem{defin}{Definition}
\newtheorem{remark}{Remark}
\newtheorem{property}{Property}
\newcommand{\la}{\frac{1}{\lambda}}
\newcommand{\sectemul}{\arabic{section}}
\renewcommand{\theequation}{\sectemul.\arabic{equation}}
\renewcommand{\thepron}{\sectemul.\arabic{pron}}
\renewcommand{\thelemma}{\sectemul.\arabic{lemma}}
\renewcommand{\thethm}{\sectemul.\arabic{thm}}
\renewcommand{\theCorol}{\sectemul.\arabic{Corol}}
\renewcommand{\theexam}{\sectemul.\arabic{exam}}
\renewcommand{\thedefin}{\sectemul.\arabic{defin}}
\renewcommand{\theremark}{\sectemul.\arabic{remark}}
\def\REF#1{\par\hangindent\parindent\indent\llap{#1\enspace}\ignorespaces}
\def\lo{\left}
\def\ro{\right}
\def\be{\begin{equation}}
\def\ee{\end{equation}}
\def\beq{\begin{eqnarray*}}
\def\eeq{\end{eqnarray*}}
\def\bea{\begin{eqnarray}}
\def\eea{\end{eqnarray}}
\def\d{\Delta_T}
\def\r{random walk}
\def\o{\overline}

\title{\large\bf On closedness under convolution roots related to an infinitely divisible distribution in the distribution class $\mathcal{L}(\gamma)$
%\thanks{Research supported by National Science Foundation of China, Grant
%No.11071182  }
}
\author{\small
Hui Xu$^{1}$~~ Yuebao Wang$^{1}$\thanks{Research supported by National Natural Science Foundation of China
(No.s 11071182, 11401415), Tian Yuan foundation (No.11426139), Natural Science Foundation of the Jiangsu
Higher Education Institutions of China (No.13KJB110025),
Postdoctoral Research Program of Jiangsu Province of China (No.
1402111C).}
\thanks{Corresponding author.
Telephone: +86 512 67422726. Fax: +86 512 65112637. E-mail:
ybwang@suda.edu.cn}~~Dongya Cheng$^{1}$~~Changjun Yu$^{2}$
\\
{\footnotesize\it 1. School of Mathematical Sciences, Soochow
University, Suzhou 215006, China}\\ {\footnotesize\it 2. School of
Sciences, Nantong University, Nantong 226019, China}}
\date{}

\maketitle

\begin{center}
{\noindent\small {\bf Abstract }}
\end{center}

{\small
%In this paper, we show that the distribution class $\mathcal{L}(\gamma)\setminus\mathcal{OS}$ for some $\gamma>0$ is not closed under convolution roots related to an infinitely divisible distribution. Precisely, we give two main conditions on L$\acute{e}$vy spectral distribution generated by L$\acute{e}$vy spectral measure of an infinitely divisible distribution, under each of which, there is an infinitely divisible distribution belonging to the class, however the corresponding spectral distribution is not. And we note that the two conditions can not be deduced from each other. For the distribution class $(\mathcal{L}(\gamma)\cap\mathcal{OS})\setminus\mathcal{S}(\gamma)$, corresponding conclusion is also proved. In order to prove the results mentioned above, we explore some of the structural properties of the class, which include the closedness of the class under the convolution and the random convolution. In addition, we study some properties of a transformation between heavy-tailed and light-tailed distributions. The transformation is a key to finding the required L$\acute{e}$vy spectral distribution. On the other hand, we also give some conditions which guarantee the closeness of the class under convolution roots. The whole research in the paper is a deep and complete discussion on the famous conjecture due to Embrechts and Goldie \cite{EG1980, EG1982} (J. Austral. Math. Soc. (Ser. A) 29, 243-256, 1980; Stoch. Process. Appl. 13, 263-278, 1982).

We consider questions related to the well-known conjecture due to Embrechts and Goldie
on the closedness of different classes of heavy- and light-tailed distributions with respect to
convolution roots.
We show that the class $\mathcal{L}(\gamma)\setminus\mathcal{OS}$ is not closed under convolution roots related to an infinitely
divisible distribution for any $\gamma\ge0$, i.e. we provide examples of infinitely divisible distributions
belonging to this class such that the corresponding Levy spectral distribution does not. We also
prove a similar statement for the class $(\mathcal{L}(\gamma)\cap\mathcal{OS})\setminus\mathcal{S}(\gamma)$.
In order to facilitate our analysis, we explore the structural properties of some of the classes
of distributions, and study some properties of the well-known transformation from a heavy-tailed distribution to a light-tailed. %one\textcolor[rgb]{1.00,0.00,0.00}{, and give a result on the closedness of the class $(\mathcal{L}(\gamma)\cap\mathcal{OS})$ under convolution roots}.
%Finally, we provide examples of rather general classes of distributions which are closed under convolution roots.
\medskip

\medskip

}

\section{Introduction and main results}
\setcounter{thm}{0}\setcounter{Corol}{0}\setcounter{lemma}{0}\setcounter{pron}{0}\setcounter{equation}{0}
\setcounter{remark}{0}\setcounter{exam}{0}\setcounter{property}{0}\setcounter{defin}{0}

In this paper, unless otherwise stated, we assume that all distributions are supported on the positive half real line.

Let $H$ be an infinitely divisible distribution
with the Laplace transform
$$
\int_0^\infty \exp\{-\lambda
y\}H(dy)=\exp\Big\{-a\lambda-\int_0^\infty(1-e^{\lambda y})\upsilon(dy)\Big\},
$$
where $a\ge0$ is a constant and the L$\acute{e}$vy spectral
measure $\upsilon$ is a Borel measure function supported on $(0,\infty)$ with the
properties $\mu=\upsilon((1,\infty))<\infty$ and $\int_0^1y\upsilon(dy)
<\infty$. Let
$$F(x)=\big(\upsilon(x)/\mu\big)\textbf{1}(x>0)=\big(\upsilon(0,x]/\mu\big)\textbf{1}(x>0),\ x\in(-\infty,\infty),$$
which is called L$\acute{e}$vy spectral distribution
generated by the measure $\upsilon$.
The distribution $H$ admits representation $H=H_1*H_2$, which is reserved
for convolution of two distributions $H_1$ and $H_2$, where
$\overline H_1(x)=1-H_1(x)=O(e^{-\beta x})$ for some $\beta>0$ and
\begin{eqnarray}\label{101}
H_2(x)=e^{-\mu}\sum_{k=0}^\infty F^{*k}(x)\mu^k/k!,\ x\in(-\infty,\infty).
\end{eqnarray}
See, for example, Feller \cite{F1971}, page 450. Here $g_1(x)=O\big(g_2(x)\big)$ means that $\limsup g_1(x)/g_2(x)<\infty$ for two positive functions $g_1$ and $g_2$.

It is well-known that one of the
important topics in the theory of infinitely divisible distribution is to discuss the closedness for certain distribution classes under
convolution roots related to an infinitely divisible distribution. {More precisely}, if an infinitely divisible distribution belongs to a certain
distribution class, and its L$\acute{e}$vy spectral distribution belongs
to the same class, then the distribution class is said to be closed under
convolution roots related to the infinitely divisible distribution.
On the contrary, if L$\acute{e}$vy spectral distribution of an infinitely divisible
distribution belongs to a certain distribution class,
and the infinitely divisible distribution  belongs to the same class,
then the distribution class is said to be closed under convolution related to the
infinitely divisible distribution.

{On the study of the topic}, we may have different {conclusions} for different distribution classes.
Thus,
%closedness under convolution roots related to a infinitely divisible distribution,
we first introduce some concepts and notations for some common distribution classes.

We say that a distribution $F$ belongs
to the distribution class $\mathcal{L}(\gamma)$ for some $\gamma\ge0$,
if for all $t$,
$$\overline {F}(x-t)\sim e^{\gamma t}\overline F(x),$$
where the notation $g_1(x)\sim g_2(x)$ means $g_1(x)/g_2(x)\to1$ for two positive functions $g_1$ and $g_2$, and all limits refer to $x$ tending to infinity.
In the above definition, if $\gamma>0$ and the distribution $F$ is lattice, then $x$ and $t$ should be restricted to values of the lattice span, see Bertoin and Doney \cite{BD1996}.
%If $F$ is a lattice point distribution, then $x$ and $t$ can only be an integer multiple of the length of the lattice point, respectively.

Further, if a distribution $F$ belongs to the class $\mathcal{L}(\gamma)$
for some $\gamma\ge0,\ \ m(F)=\int_0^\infty e^{\gamma y}F(dy)<\infty$ and
$$\overline {F^{*2}}(x)\sim2m(F)\overline F(x),$$
then we say that the distribution $F$ belongs to the distribution class $\mathcal{S}(\gamma)$.

In particular, the classes $\mathcal{L}=\mathcal{L}(0)$ and $\mathcal{S}=\mathcal{S}(0)$ are called the long-tailed distribution class and the subexponential distribution class, respectively.
%We note that, in the definition of class $\mathcal{L}(\gamma)$, when $\gamma=0$, $\overline {F}(x-t)\sim \overline F(x)$ for all $t \neq 0$ can be substituted by $\overline {F}(x-t)\sim \overline F(x)$ for some $t \neq 0$. Combined with Lemma 2.2 of \cite{YWC2010} and (\ref{303}) of this paper, however, we know that the similar conclusion does not hold when $\gamma>0$.
It should be noted that the requirement $F\in\mathcal{L}$ is not needed in the definition of the class $\mathcal{S}$.
The class $\mathcal{L}$ is a heavy-tailed distribution subclass, while $\mathcal{L}(\gamma)$ for $\gamma>0$ is
a light-tailed distribution subclass.
We denote the heavy-tailed distribution class by $\mathcal{K}$, and the light-tailed distribution class by $\mathcal{K}^c$.

The classes $\mathcal{L}(\gamma)$ and  $\mathcal{S}(\gamma)$
were introduced by Chistyakov \cite{C1964} for $\gamma=0$ and Chover et al. \cite{CNW1973a,CNW1973b} for $\gamma>0$, respectively. And the class $\cup_{\gamma\ge0}\mathcal{L}(\gamma)$ is properly contained
in the following distribution class introduced by Shimura and Watanabe \cite{SW2005-1}.

We say that a distribution $F$ belongs
to the generalized long-tailed distribution class, denoted by $F\in\mathcal{OL}$,
if for all (or, equivalently, for some) $t \neq 0$,
$$C(F,t)=\limsup\overline {F}(x-t)/\overline F(x)<\infty.$$

Correspondingly, the class $\cup_{\gamma\ge0}\mathcal{S}(\gamma)$ is properly contained
in the following distribution class introduced by Kl\"{u}ppelberg \cite{K1990}.
We say that a distribution $F$ belongs
to the generalized subexponetial distribution class, denoted by $F\in\mathcal{OS}$, if
$$C^*(F)=\limsup\overline {F^{*2}}(x)/\overline F(x)<\infty,$$
see also \cite{SW2005-1} for more details about the class.

%The following distribution class is a subclass of the subexponential distribution class.

%We say that the distribution $F$ belongs to the strong subexponential distribution class introduced by Kl\"{u}ppelberg \cite{K1988}, denoted by $F\in\mathcal{S}^*$, if $0<\mu(F)=\int_0^\infty \overline{F}(y)dy<\infty$ and $$\int_0^x \overline{F}(x-y)\overline{F}(y)dy\sim2\mu(F)\overline{F}(x).$$

%Wang et al. \cite{WXCY2014} extended the class $\mathcal{S}^*$ to the following distribution class.

%We say that the distribution $F$ belongs to the generalized strong subexponential distribution class, denoted by $F\in\mathcal{OS}^*$, if $$C^\otimes(F)=\limsup\int_0^x \overline{F}(x-y)\overline{F}(y)dy/\overline F(x)<\infty.$$

Apart from the distribution classes, {the study of the above topic} also depends on the relationship between $H_1$ and $H_2$, which were introduced before. In this paper, we will give a wider range of choices for them. For example, we can consider the case that, for some $C\in[0,\infty)$,
\begin{eqnarray}\label{102}
\lim\overline{H_1}(x)/\overline{H_2}(x)=C.
\end{eqnarray}
More generally, we consider the two distributions $H_1$ and $H_2$ such that
\begin{eqnarray}\label{103}
\mid\overline{H_1}(x-t)-e^{\gamma t}\overline{H_1}(x)\mid = o(\overline{H_2}(x)),
\end{eqnarray}
{where $g_1(x)=o(g_2(x))$ means that $g_1(x)/g_2(x)\to0$ for two positive functions $g_1$ and $g_2$}. When $H_2\in\mathcal{L}(\gamma)$ and (\ref{102}) is satisfied, then (\ref{103}) holds. On the contrary, if (\ref{103}) is satisfied, then (\ref{102}) does not necessarily hold, see Subsection 6.1 {for details}.

For the class $\mathcal{S}(\gamma)$ with $\gamma\ge0$, the closedness
under convolution roots related to an infinitely
divisible distribution has been proved,
%with condition that the spectral distribution$ F\in\mathcal{L}(\gamma)$ for some $\gamma>0$,
see Embrechts et al. \cite{EGV1979} for $\gamma=0$ and {Sgibnev \cite{S1990} or Pakes \cite{P2004}} for $\gamma>0$.
%Embrechts and Goldie \cite{EG1982} for $\gamma>0$ with condition that the L$\acute{e}$vy spectra distribution belongs to the class $\mathcal{L}(\gamma)$.
Thus, an interesting problem is raised naturally:
$$If\ H\in\mathcal{L}(\gamma),\ does\ this\ imply\ that\ F\in\mathcal{L}(\gamma)?$$
In order to clarify this issue, we might firstly answer the following problem:
$$If\ H_2\in\mathcal{L}(\gamma),\ does\ this\ imply\ that\ F\in\mathcal{L}(\gamma)?$$
For the class $\mathcal{L}(\gamma)$, the latter problem on closedness
under random convolution roots is a natural extension of the famous Embrechts-Goldie's
conjecture on closedness of the class $\mathcal{L}(\gamma)$ under convolution roots, see \cite{EG1980, EG1982}:
$$If\ F^{*k}\in\mathcal{L}(\gamma)\ for\ some\ (even\ for\ all)\ k\ge2,\ then\ F\in\mathcal{L}(\gamma).$$
Therefore, we call the above-mentioned two problems the generalized Embrechts-Goldie's problems.

So far, Shimura and Watanabe \cite{SW2005-2} {have given a negative answer to} the Embrechts-Goldie's conjecture in the case that $\gamma\ge0$ and $k=2$. Furthermore, Watanabe \cite{W2015} showed that the class $\mathcal{S}(\gamma)$ with $\gamma>0$ is not closed under convolution roots, and got some related important results. However, they did not discuss  the generalized Embrechts-Goldie's problems. Recently, in the case that $\gamma=0$ and $k\ge2$, Xu et al. \cite{XFW2015} gave a negative answer to the Embrechts-Goldie's conjecture and the above-mentioned two problems for the classes $\mathcal{L}\setminus\mathcal{OS}$ and $\mathcal{L}\cap\mathcal{OS}\setminus\mathcal{S}$. Therefore, the generalized Embrechts-Goldie's problems in the case that $\gamma>0$ attracted our attention. %At first glance, these problems in the case that $\gamma>0$ is similar to problems in the case that $\gamma=0$, however,

We soon find that the method in the case $\gamma=0$ can not be used in the case $\gamma>0$ directly. Thus, we need to find a new method with more technical details to provide negative answers to the generalized Embrechts-Goldie's problems for the classes $\mathcal{L}(\gamma)\setminus\mathcal{OS}$ and $(\mathcal{L}(\gamma)\cap\mathcal{OS})\setminus\mathcal{S}(\gamma)$, respectively.
%In addition, this paper shows that the classes and $\mathcal{OL}\cap\mathcal{K}^{c}$ are not closed under convolution roots related to an infinitely divisible distribution too.

%For the sake of completeness, the following result is discussed in the case that $\gamma\ge0$, among them, the result for the case that $\gamma=0$ with different methods belongs to Theorem 1.2 (3) in \cite{XFW2015}.

%\begin{thm}\label{thm101} For any $\gamma\geq 0$, there is an infinitely divisible distribution $H$ such that $H\in\cal{L}(\gamma)$, while its L$\acute{e}$vy spectra distribution $F\in\mathcal{OL}\setminus\mathcal{L}(\gamma)$, even, $F\notin\cal{OL}$. \end{thm}

\begin{thm}\label{thm101}
For any $\gamma>0$,
%there are four {families} of L$\acute{e}$vy spectral distribution $F$ of an infinitely divisible distribution $H$ satisfying condition
assume that the L$\acute{e}$vy spectral distribution $F$ of an infinitely divisible distribution $H$ satisfies condition
\begin{eqnarray}\label{104}
\liminf\overline F(x-t)/\overline F(x)\ge e^{\gamma t}\ for\ all\ t>0,
\end{eqnarray}
or
\begin{eqnarray}\label{1040}
\overline{F}(x)=o(\overline{F^{*2}}(x)).
\end{eqnarray}
Among them, there is a L$\acute{e}$vy spectral distribution $F$ such that $F\in\mathcal{OL}\setminus\mathcal{L}(\gamma)$ with $m(F)=\infty$, or, $F\notin\cal{OL}$, while $H$ and $F^{*k}$ for all $k\ge2$ belong to the class $\cal{L}(\gamma)\setminus\mathcal{OS}$;
%In the remaining family,
and there is an another L$\acute{e}$vy spectral distribution $F$ such that $F\in\cal{OL}\setminus\big(\mathcal{L}(\gamma)\cup\mathcal{OS}\big)$ with $m(F)<\infty$, while $H$ for all $k\ge2$ and $F^{*k}$ belong to the class $\big(\mathcal{L}(\gamma )\cap\mathcal{OS}\big)\setminus\mathcal{S}(\gamma)$, as well as
\begin{eqnarray}\label{10400}
&&e^{-\mu}\sum_{m=1}^\infty \Big(\sum_{k=2m}^{2(m+1)-1}\mu^k/k!\Big) m(F^{*2})^{m-1}m\le\liminf\overline{H}(x)/\overline{F^{*2}}(x)\le\limsup\overline{H}(x)/\overline{F^{*2}}(x)\nonumber\\
&\le&e^{-\mu}\sum_{m=1}^\infty\Big(\sum_{k=2(m-1)+1}^{2m}\mu^k/k! \Big)\sum_{i=0}^{m-1}\big(m(F^{*2})\big)^i\big(C^*(F^{*2})-m(F^{*2})\big)^{m-1-i}.
\end{eqnarray}
\end{thm}

\begin{remark}\label{rem103}
%From {Theorem \ref{thm101}} and Remark \ref{101}, we know that, the conditions (\ref{104}) and (\ref{1040}), more generally, the conditions (\ref{-1}) and (\ref{400}), are two main conditions in the present paper. Among them,
The condition (\ref{104}) has been used in Lemma 7 and Theorem 7 of Foss and Korshunov \cite{FK2007}, Wang et al. \cite{WXC2016}, and so on.
%In addition, the condition (\ref{104}) for all $t>0$ can not be replaced by the condition (\ref{104}) for some $t>0$, see Remark 5.3 of Yu et al. \cite{YWC2010}.

Some distributions which do not satisfy the conditions (\ref{104}) and (\ref{1040}) have been found in Example 1 in \cite{FK2007}, Proposition 3.2 and Remark 4.1 in Chen et al. \cite{CYW2013}, Theorem 1.1 in \cite{W2015}, and so on. On the other hand, many distributions satisfy condition (\ref{104}) or condition (\ref{1040}). Clearly, condition (\ref{104}) is satisfied when $\gamma=0$ or $F\in\mathcal{L}(\gamma)$ for all $\gamma>0$. In addition, when $\gamma>0$ and $F\notin\mathcal{L}(\gamma)$ for all $\gamma>0$, there are some distributions which simultaneously satisfy conditions (\ref{104}) and (\ref{1040}), for example, the distribution $F$ in classes $\mathcal{F}_i(\gamma)$ for some $\gamma>0,\ i=1,2$, see Definition \ref{defin301}, Propositions \ref{pron301} and \ref{pron302} below. However, the distributions in \cite {SW2005-2} and in the class $\mathcal{F}_4(\gamma)$ with some $\gamma>0$ below satisfy (\ref{1040}) but do not satisfy (\ref{104}); and vice versa, the distributions in the class $\mathcal{F}_3(\gamma)$ with some $\gamma>0$ satisfy (\ref{104}) but do not satisfy (\ref{1040}), see Propositions \ref{pron303} and  \ref{pron304} below.
In other words, the conditions (\ref{104}) and (\ref{1040}) can not be deduced from each other.
%in other words, these two conditions have their own independent value of existence.
\end{remark}

Theorem \ref{thm101} %andProposition \ref{pron101} with Corollary \ref{Corol102} are
is proved in Section 5. To this end, some structural properties of the class $\cal{L}(\gamma)$ for some $\gamma\ge0$, as preliminary results, are respectively given in Section 2 and Section 3, which include the closedness of the class under the convolution and the convolution roots. In Section 4, a transformation from heavy-tailed distributions to light-tailed distributions is studied. Finally, some remarks on the above-mentioned results, as well as the relevant conditions, and a local version on the study of Embrechts-Goldie¡¯s conjecture, are given in Section 6.

\section{On the closedness under the convolution}
\setcounter{thm}{0}\setcounter{Corol}{0}\setcounter{lemma}{0}\setcounter{pron}{0}\setcounter{equation}{0}
\setcounter{remark}{0}\setcounter{exam}{0}\setcounter{property}{0}\setcounter{defin}{0}

\date{}
Before presenting the main results in this section, we first provide the following lemma, which together with Remark \ref{rem201} below, is also the key to proving the main results in this paper and has its own independent interest.

\begin{lemma}\label{lemma202}
Let $ F_1$ and $F_2$ be two distributions such that $F_1*F_2 \in \mathcal{L}(\gamma)$ for some $\gamma\ge0$.
Further, assume that
\begin{eqnarray}\label{-1}
\liminf\overline F_i(x-t)/\overline F_i(x)\ge e^{\gamma t}\ \ for\ all\ t>0,\ i=1,2,
\end{eqnarray}
then
\begin{eqnarray}\label{0}
\mid\overline F_i(x-t)-e^{\gamma t}\overline F_i(x)\mid=o(\overline {F_1*F_2}(x))\ \ for\ all\ t>0,\ i=1,2.
\end{eqnarray}
Particularly, let $F$ be a distribution such that $F^{*n}\in\cal{L}(\gamma)$ for some integer $n\ge2$ and some $\gamma\ge0$. Further, if
\begin{eqnarray}\label{00}
\liminf\overline {F^{*i}}(x-t)/\overline {F^{*i}}(x)\ge e^{\gamma t}\ for\ \ all\ t>0,\ 1\le i\le n-1,
\end{eqnarray}
then
\begin{eqnarray}\label{000}
\mid\overline {F^{*i}}(x-t)-e^{\gamma t}\overline {F^{*i}}(x)\mid=o(\overline {F^{*n}}(x))\ for\ \ all\ t>0,\ 1\le i\le n-1.
\end{eqnarray}
\end{lemma}

\begin{remark}\label{rem201}
In the case of $\gamma=0$, the condition (\ref{00}) holds automatically, the conclusion (\ref{0}) and some of the following results are due to Theorems 1.1 and 1.2 {in} \cite{XFW2015}. However, as we point out, our method here is different from that in \cite{XFW2015}.
%It should be noted that the method used in case $\gamma>0$ is different from the method in case $\gamma=0$.

In addition, under the condition $F_1*F_2 \in \mathcal{L}(\gamma)$ for some $\gamma>0$, if the condition (\ref{-1}) is substituted by the following condition
\begin{eqnarray}\label{400}
\overline {F_i}(x)=o(\overline{F_1*F_2}(x)),\ i=1,2,
\end{eqnarray}
then the conclusion (\ref{0}) still holds, obviously.

%In particular, if $F_1=F_2=F$, then (\ref{-1}) and (\ref{400}) are reduced to (\ref{104}) and (\ref{1040}) \textcolor[rgb]{1.00,0.00,0.00}{in Section 3, among which, the condition (\ref{104}) for all $t>0$ have been used in Lemma 7 and Theorem 7 of Foss and Korshunov \cite{FK2007}, and so on.}
%\begin{eqnarray}\label{1040} \overline{F}(x)=o(\overline{F^{*2}}(x)), \end{eqnarray}
\end{remark}
\proof Firstly, we prove that for $1\le i\neq j\le2,$

\begin{eqnarray}\label{5}
\lim_{v\to\infty}\limsup\int_{[0,v]}\big(\overline F_i(x-t-y)-e^{\gamma t}\overline F_i(x-y)\big)F_j(dy)/\overline{F_1*F_2}(x)=0
\end{eqnarray}
and
\begin{eqnarray}\label{50}
\lim_{v\to\infty}\liminf\int_{[0,v]}\big(\overline F_i(x-t-y)-e^{\gamma t}\overline F_i(x-y)\big)F_j(dy)/\overline{F_1*F_2}(x)=0.
\end{eqnarray}

We point out that, when $\gamma = 0$, the conclusions (\ref{5}) and (\ref{50})
are implied by Theorem 1.1 (1a) of \cite{XFW2015}.
%In fact, for any $t,v>0$ and $1\le i\neq j\le2$, from the above result, we have
%$$\lim\int_{0}^{v}\Big(\overline F_i(x-t-y)-\overline F_i(x-y)\Big)F_j(dy)/\overline{F_1*F_2}(x)=0.$$
Thus we only need to prove (\ref{5}) and (\ref{50}) in the case that $i=1$ and $\gamma>0$. %\textcolor[rgb]{1.00,0.00,0.00}{And without loss of generality, there we assume that $F_i$ is absolutely continuous, $i=1,2$.}

For any $x>v$, using integration by parts, we have
\begin{eqnarray}\label{1}
&&\overline {F_1*F_2}(x-t)=\overline F_2(x-t)+\Big(\int_{[0,x-t-v]}+\int_{(x-t-v,x-t]}\Big)\overline F_1(x-t-y)F_2(dy)\nonumber\\
&=&e^{\gamma t}\overline {F_1*F_2}(x)+\overline F_2(x-t)+\int_{[0,x-t-v]}\big(\overline F_1(x-t-y)-e^{\gamma t}\overline F_1(x-y)\big)F_2(dy)\nonumber\\
& &+\int_{(x-t-v,x-t]}\overline F_1(x-t-y)F_2(dy)-e^{\gamma t}\overline F_2(x)-e^{\gamma t}\int_{(x-t-v,x]}\overline F_1(x-y)F_2(dy)\nonumber\\
&=&e^{\gamma t}\overline {F_1*F_2}(x)+\int_{[0,v]}\big(\overline F_2(x-t-y)-e^{\gamma t}\overline F_2(x-y)\big)F_1(dy)-\int_{(v,v+t]}e^{\gamma t}\overline F_2(x-y)F_1(dy)\nonumber\\
& &+\int_{[0,x-t-v]}\big(\overline F_1(x-t-y)-e^{\gamma t}\overline F_1(x-y)\big)F_2(dy)+
\big(\overline {F_1}(v)-e^{\gamma t}\overline {F_1}(v+t)\big)\overline F_2(x-t-v)\nonumber\\
&\geq&e^{\gamma t}\overline {F_1*F_2}(x)+\int_{[0,v]}\big(\overline F_2(x-t-y)-e^{\gamma t}\overline F_2(x-y)\big)F_1(dy)-(e^{\gamma t}-1)\overline F_2(x-v-t)\overline F_1(v)\nonumber\\
& &+\int_{[0,x-t-v]}\big(\overline F_1(x-t-y)-e^{\gamma t}\overline F_1(x-y)\big)F_2(dy).
\end{eqnarray}

In the following, the notation $g_1(x)\gtrsim\ (or\ \lesssim) g_2(x)$ means that $\limsup g_2(x)/g_1(x)$ (or $g_1(x)/g_2(x))\le1$ for two positive functions $g_1$ and $g_2$ supported on $[0,\infty)$. Now, we deal with\\
$(e^{\gamma t}-1)\overline F_2(x-v-t)\overline F_1(v)$ in (\ref{1}).
By $F_1*F_2\in\mathcal{L}(\gamma)$, (\ref{-1}) and Fatou$^{,}$s lemma,

\begin{eqnarray*}
&&\overline {F_1*F_2}(x)/(\overline F_2(x-v-t)\overline {F_1}(v))\sim\overline {F_1*F_2}(x-v-t)e^{-\gamma(v+t)}/(\overline F_2(x-v-t)\overline {F_1}(v))\nonumber\\
&\gtrsim&\Big(\sum_{k=1}^{[v/t]}\int_{(v-kt,v-(k-1)t]}e^{\gamma y}F_1(dy)\Big)/(e^{\gamma (v+t)}\overline {F_1}(v))\nonumber\\
&\geq&\Big(\sum_{k=1}^{[v/t]}e^{\gamma (v -kt)}\big(\overline F_1(v-k t)-\overline F_1(v-k t+t)\big)\Big)/(e^{\gamma (v+t)}\overline {F_1}(v))=a(v),
\end{eqnarray*}
where $[c]$ is the integer part of the number $c$. Then, still by (\ref{-1}) and Fatou$^{,}$s lemma,
\begin{eqnarray*}
\liminf_{v\to\infty}a(v)&=&\liminf_{v\to\infty}\sum_{k=1}^{[v/t]}e^{-\gamma (k+1)t}\big(\overline F_1(v-k t)/\overline F_1(v-k t+t)-1\big)\big(\overline F_1(v-kt+t)/\overline {F_1}(v)\big)\nonumber\\
&\ge&\sum_{k=1}^\infty e^{-\gamma t}(1-e^{-\gamma t})=\infty.
\end{eqnarray*}
Thus,
\begin{eqnarray}\label{2}
\lim_{v\to\infty}\limsup(e^{\gamma t}-1)\overline F_2(x-v-t)\overline F_1(v)/\overline {F_1*F_2}(x)=0.
\end{eqnarray}

According to the condition (\ref{-1}), for any $\epsilon>0$, we know that there is a constant $v_0>0$ such that for $v\ge v_0$ and $x\ge2v+t$,

\begin{eqnarray}\label{3}
\int_{[0,x-v-t]}\big(\overline F_1(x-t-y)-e^{\gamma t}\overline F_1(x-y)\big)F_2(dy)
&\wedge&\int_{[0,v]}\big(\overline F_1(x-t-y)-e^{\gamma t}\overline F_1(x-y)\big)F_2(dy)\nonumber\\
\ge-\epsilon\int_{[0,x]}\overline F_1(x-y)F_2(dy)\ge-\epsilon\overline {F_1*F_2}(x),
\end{eqnarray}
where $a\wedge b=\min\{a,b\}$. On the other hand, by $F_1*F_2\in\cal{L}(\gamma)$ and (\ref{1})-(\ref{3}), we have
\begin{eqnarray}\label{4}
\lim_{v\to\infty}\limsup\int_{[0,v]}(\overline F_2(x-t-y)-e^{\gamma t}\overline {F_2} (x-y))/\overline{F_1*F_2}(x)F_1(dy)\le\epsilon.
\end{eqnarray}

From (\ref{3}), (\ref{4}) and the arbitrariness of $\epsilon$, we immediately obtain (\ref{5}) and (\ref{50}).

Secondly, we only need to prove (\ref{0}) in the case that $i=1$ and $\gamma>0$.
%At this time, we have
%\begin{eqnarray*}
%\liminf(\overline {F_1}(x-t)-e^{\gamma t}\overline F_1(x))/\overline {F}(x)\ge0.
%\end{eqnarray*}
To this end, we prove
\begin{eqnarray}\label{511}
\liminf\big(\overline {F_1}(x-t)-e^{\gamma t}\overline F_1(x)\big)/\overline {F_1*F_2}(x)=0.
\end{eqnarray}
Otherwise, there exists a constant $C>0$ such that
\begin{eqnarray}\label{52}
\liminf\big(\overline {F_1}(x-t)-e^{\gamma t}\overline F_1x)\big)/\overline {F_1*F_2}(x)\ge C.
\end{eqnarray}
Then (\ref{50}), (\ref{52}), Fatou's Lemma and $F_1*F_2\in\cal{L}(\gamma)$ lead to the following contradiction:
\begin{eqnarray*}
&&0=\lim_{v\to\infty}\liminf\int_{[0,v]}\big(\overline F_1(x-t-y)-e^{\gamma t}\overline {F_1}(x-y)\big)F_2(dy)/\overline {F_1*F_2}(x)\nonumber\\
&\ge&\lim_{v\to\infty}\int_{[0,v]}\liminf\frac{\overline F_1(x-t-y)-e^{\gamma t}\overline F_1(x-y)}{\overline {F_1*F_2}(x-y)}\liminf\frac{\overline {F_1*F_2}(x-y)}{\overline {F_1*F_2}(x)}F_2(dy)\nonumber\\
&\ge&C\int_{[0,\infty)}e^{\gamma y}F_2(dy)>0.
\end{eqnarray*}
Therefore (\ref{511}) holds.

Now, we prove
\begin{eqnarray}\label{54}
\limsup\big(\overline {F_1}(x-t)-e^{\gamma t}\overline {F_1}(x)\big)/\overline {F_1*F_2}(x)\le0.
\end{eqnarray}
By (\ref{5}), Fatou's Lemma, (\ref{511}) and $F_1*F_2\in\mathcal{L}(\gamma)$, we have

\begin{eqnarray}\label{10}
&&0=\lim_{v\to \infty}\limsup\Big(\int_{[0,t)}+\int_{[t,2t)}+\int_{[2t,v]}\Big)\big(\overline {F_1}(x-2t-y)-e^{2\gamma t}\overline {F_1}(x-y)\big)/\overline {F_1*F_2}(x)F_2(dy)\nonumber\\
&\ge&\limsup\int_{[t,2t)}\big(\overline {F_1}(x-2t-y)-e^{\gamma y}\overline {F_1}(x-2t)\big)F_2(dy)/\overline {F_1*F_2}(x)\nonumber\\
&&\ \ \ \ \ \ +\int_{[t,2t)}\liminf e^{\gamma y}\big(\overline {F_1}((x-y)-(2t-y))-e^{\gamma (2t-y)}\overline {F_1}(x-y)\big)F_2(dy)/\overline {F_1*F_2}(x)\nonumber\\
&&\ \ \ \ \ \ +\lim_{v\to \infty}\Big(\int_{[0,t)}+\int_{[2t,v]}\Big)e^{\gamma y}\liminf\frac{\overline {F_1}\big((x-y)-2t)-e^{2\gamma t}\overline {F_1}(x-y)}{\overline {F_1*F_2}(x-y)}F_2(dy)\nonumber\\
&=&\limsup\int_{[t,2t)}\big(\overline {F_1}(x-2t-y)-e^{\gamma y}\overline {F_1}(x-2t)\big)F_2(dy)/\overline {F_1*F_2}(x)\nonumber\\
&\geq&\int_{[t,2t)}\liminf\big(\overline {F_1}((x-3t)-(y-t))-e^{\gamma (y-t)}\overline {F_1}(x-3t)\big)F_2(dy)/\overline {F_1*F_2}(x-3t)\nonumber\\
&&\ \ \ \ \ \ +\limsup\int_{[t,2t)} e^{\gamma (y-t)}\big(\overline {F_1}((x-2t)-t)-e^{\gamma t}\overline {F_1}(x-2t)\big)F_2(dy)/\overline {F_1*F_2}(x-2t)\nonumber\\
&=&e^{\gamma t}\int_{[t,2t)} e^{\gamma y}F_2(dy)\limsup\big(\overline {F_1}((x-2t)-t)-e^{\gamma t}\overline {F_1}(x-2t)\big)/\overline {F_1*F_2}(x-2t).
\end{eqnarray}
Thus (\ref{54}) follows from (\ref{10}).

Combining (\ref{54}) with (\ref{511}) yields (\ref{0}) with $i=1$.\hfill$\Box$
\vspace{0.2cm}

In the following, we give a result on closedness of the class $\mathcal{L}(\gamma)$ under convolution.
\begin{thm}\label{thm201}
For any integer $n\ge2$, let $F_i,\ 1\le i\le n$, be $n$ distributions such that $F_n\in\mathcal{L}(\gamma)$ for some $\gamma\ge0$. Further, assume that condition (\ref{-1}) and
\begin{eqnarray}\label{130}
\mid\overline{F_i}(x-t)-e^{\gamma t}\overline{F_i}(x)\mid= o\big(\overline{F_n}(x)\big)\ \ for\ all\ t>0\ and\ i=1,\cdot\cdot\cdot,n-1,
\end{eqnarray}
are satisfied.
%Or, assume that
%\begin{eqnarray}\label{131}
%\overline{F_i}(x) = o(\overline{F_n}(x)),\ i=1,\cdot\cdot\cdot,n-1.
%\end{eqnarray}
Then $F_1\ast\cdot\cdot\cdot\ast F_n\in\mathcal{L}(\gamma)$.
\end{thm}

\begin{remark}\label{rem203}
For $n=2$ and some $\gamma\ge0$, Theorem 3 of Embrechts and Goldie \cite{EG1980} asserted the following result on the closedness of the class $\cal{L}(\gamma)$ under convolution:
assume $F_2\in\mathcal{L}(\gamma)$, then $F=F_1*F_2\in \mathcal{L}(\gamma)$
if either (a) $F_1\in\mathcal{L}(\gamma)$ or (b) $\overline{F_1}(x)=o\big(\overline{F_2}(x)\big)$.
{For the latter}, more generally, if $F_n\in\mathcal{L}(\gamma)$ and
\begin{eqnarray}\label{1300}
\overline{F_i}(x) = o\big(\overline{F_n}(x)\big),\ 1\le i\le n-1,
\end{eqnarray}
then $F_1\ast\cdot\cdot\cdot\ast F_n\in\mathcal{L}(\gamma)$.
Clearly, (\ref{130}) may be implied by (\ref{1300}) and $F_n\in\mathcal{L}(\gamma)$, but the inverse implication dose not always hold, see Subsection 6.1 for the details, so
%condition (\ref{130}) is a rather loose condition, see Remark \ref{rem103} and Subsection 6.1 for the details. However, as discussed in Remark \ref{rem103}, the conditions in (\ref{-1}) with $1\le i\le n-1$ are different from the condition (\ref{1300}). Therefore,  as Theorem 1.1 (1b) or Remark 2.1 {in} \cite{XFW2015} for the case $\gamma=0$,
the theorem is a slight extension of Theorem 3 in \cite{EG1980} for the case that $\gamma>0$.
\end{remark}
\proof It is sufficient to prove the result for $n = 2$ in view of the
induction method. We take $v$ large enough and $x>2v+1$. By $F_2\in\cal{L}(\gamma)$, (\ref{1}) and integration by parts, we have
\begin{eqnarray}\label{13}
&&\overline {F_1*F_2}(x-t)=e^{\gamma t}\overline {F_1*F_2}(x)+\int_{[0,v]}\big(\overline F_1(x-t-y)-e^{\gamma t}\overline F_1(x-y)\big)F_2(dy)\nonumber\\
&&\ \ \ \ \ \ +\big(\overline {F_2}(v)-e^{\gamma t}\overline {F_2}(v+t)\big)\overline F_1(x-t-v)-\int_{(v,v+t]}e^{\gamma t}\overline F_1(x-y)F_2(dy)\nonumber\\
&\leq&e^{\gamma t}\overline {F_1*F_2}(x)+\int_{[0,v]}\big(\overline F_1(x-t-y)-e^{\gamma t}\overline F_1(x-y)\big)F_2(dy)\nonumber\\
& &\ \ \ \ \ \ +\big(\overline {F_2}(v)-e^{\gamma t}\overline {F_2}(v+t)\big)\overline F_1(x-t-v)-e^{\gamma t}\overline F_1(x-v)\big(\overline {F_2}(v)-\overline {F_2}(v+t)\big)\nonumber\\
&\leq&e^{\gamma t}\overline {F_1*F_2}(x)+\int_{[0,v]}\big(\overline F_1(x-t-y)-e^{\gamma t}\overline F_1(x-y)\big)F_2(dy)\nonumber\\
&&\ \ \ \ \ \ +3e^{\gamma t}\overline {F_2}(v+t)\overline F_1(x-t-v).
\end{eqnarray}

Similar to the proof of (\ref{2}), we know by (\ref{-1}) and $F_2\in\cal{L}(\gamma)$ that

\begin{eqnarray}\label{14}
&&\overline {F_1*F_2}(x)/\big(\overline {F_2}(v+t)\overline F_1(x-t-v)\big)\ge\int_{(v+t,2v]}\overline {F_1}(x-y)F_2(dy)/\big(\overline {F_2}(v+t)\overline F_1(x-t-v)\big)\nonumber\\
&\gtrsim&\int_{(v+t,2v]}e^{\gamma y}F_2(dy)/\big(e^{\gamma (v+t)}\overline {F_2}(v+t)\big)\nonumber\\
&\geq&\sum_{k=1}^{[vt^{-1}]}e^{\gamma (v +kt)}\big(\overline {F_2}(v+k t)-\overline {F_2}(v+k t+t)\big)/\big(e^{\gamma (v+t)}\overline {F_2}(v+t)\big)\nonumber\\
&\sim&e^{-\gamma t}\sum_{k=1}^{[vt^{-1}]}e^{\gamma kt}(e^{-\gamma(k -1)t}-e^{-\gamma k t})\nonumber\\
&=&[vt^{-1}](1-e^{-\gamma t})\to\infty,~~~v\to\infty.
\end{eqnarray}

According to (\ref{-1}), for any $\epsilon>0$, there is a constant $v_0>0$ such that for $v\ge v_0$ and $x\ge2v+t$,
\begin{eqnarray}\label{15}
\int_{[0,v]}\big(\overline F_1(x-t-y)-e^{\gamma t}\overline F_1(x-y)\big)F_2(dy)&\ge&-\epsilon\int_{[0,v]}\overline {F_1}(x-y)F_2(dy)\nonumber\\
&\ge&-\epsilon\overline {F_1*F_2}(x).
\end{eqnarray}
And by (\ref{130}),

\begin{eqnarray}\label{16}
&&\int_{[0,v]}\big(\overline F_1(x-t-y)-e^{\gamma t}\overline F_1(x-y)\big)F_2(dy)
=\int_{[0,v]}\big(\overline F_1(x-t-y)-e^{\gamma y}\overline F_1(x-t)\big)F_2(dy)\nonumber\\
&&\ \ \ \ \ \ +\int_{[0,t)}e^{\gamma y}\big(\overline F_1((x-y)-(t-y))-e^{\gamma (t-y)}\overline F_1(x-y)\big)F_2(dy)\nonumber\\
&&\ \ \ \ \ \ -\int_{[t,v]}e^{\gamma t}\big(\overline F_1((x-t)-(y-t))-e^{\gamma (y-t)}\overline F_1(x-t)\big)F_2(dy)\nonumber\\
&\le&\epsilon\int_{[0,v]}\overline {F_2}(x-t)F_2(dy)
+\epsilon\int_{[0,t)}e^{\gamma y}\overline {F_2}(x-y)F_2(dy)
+\epsilon\int_{[t,v]}e^{\gamma t}\overline {F_2}(x-t)F_2(dy)\nonumber\\
&\lesssim&\epsilon\big(e^{\gamma t}+2e^{2\gamma t}\big)\overline {F_1*F_2}(x).
\end{eqnarray}
Combining (\ref{13}), (\ref{14}), (\ref{15}) and (\ref{16}), we can get $F_1*F_2\in\mathcal{L}(\gamma)$.\hfill$\Box$
\vspace{0.2cm}

Based on Theorem \ref{thm201}, we obtain the following result.
\begin{Corol}\label{Corol201}
For any integer $n\ge2$ and $\gamma\ge0$, suppose that $F^{*n} \in \mathcal{L}(\gamma)$ and either of the
the following two cases is true:\\
Case (1). Condition (\ref{104})
%\begin{eqnarray}\label{17}
%\liminf\overline F(x-t)/\overline F(x)\ge e^{\gamma t}
%\end{eqnarray}
and the condition that
\begin{eqnarray}\label{18}
\mid\overline{F}(x-t)-e^{\gamma t}\overline{F}(x)\mid=o\big(\overline{F^{*n}}(x)\big)\ for\ all\ t>0
\end{eqnarray}
hold;\\
Case (2). The condition (\ref{1040}) is satisfied.\\
%\begin{eqnarray}\label{190}
%\overline {F}(x)=o(\overline {F^{*2}}(x)).
%\end{eqnarray}
Then $F^{*k}\in\mathcal{L}(\gamma)$ for all integers $k\ge n+1$.
\end{Corol}

\begin{remark}\label{rem204}
According to $F^{*n}\in\mathcal{L}(\gamma)$, (\ref{104}) and Lemma \ref{lemma202}, we know that the condition (\ref{18}) can be implied by the following condition
\begin{eqnarray}\label{19}
\liminf\overline {F^{*(n-1)}}(x-t)/\overline {F^{*(n-1)}}(x)\ge e^{\gamma t}\ for\ all\ t>0.
\end{eqnarray}
Specifically, when $n=2$, (\ref{19}) is just the same to (\ref{104}). In addition, if condition (\ref{104}) is replaced by condition (\ref{00}), then condition (\ref{18}) can be canceled according to Lemma \ref{lemma202}.
\end{remark}

\proof According to the mathematical induction method, we only need to prove that $F^{*(n+1)}\in\mathcal{L}(\gamma)$.

In Case (1), by $F^{*n}\in\mathcal{L}(\gamma)$, we know that for all $t>0$,
$$\liminf\overline {F^{*n}}(x-t)/\overline {F^{*n}}(x)\ge e^{\gamma t}$$
and
$$\mid\overline{F^{*n}}(x-t)-e^{\gamma t}\overline {F^{*n}}(x)\mid=o(\overline{F^{*n}}(x)).$$
Further, by (\ref{104}), (\ref{18}) and Theorem \ref{thm201}, we know that
$F^{*(n+1)}=F*F^{*n}\in \mathcal{L}(\gamma)$.

In Case (2), we first prove the following result:
if condition (\ref{1040}) is satisfied, then
\begin{eqnarray}\label{191}
\overline {F^{*k}}(x)=o\big(\overline {F^{*(k+1)}}(x)\big),\ k\ge1.
\end{eqnarray}
In fact, according to (\ref{1040}) and the induction hypothesis, for any $\epsilon>0$, there is a constant $x_k>0$ such that  $\overline{F}(x)\leq\varepsilon\overline{F^{*2}}(x)$ and $\overline{F^{*(k-1)}}(x)\leq\varepsilon\overline{F^{*k}}(x)$ hold for $x\ge x_k$. Thus
 using the method of integration by parts, we know that

\begin{eqnarray*}
\overline{F^{*k}}(x)&=&\overline {F}(x)+\lo(\int_{[0,x-x_k]}+\int_{(x-x_k,x]}\ro)\overline {F^{*(k-1)}}(x-y)F(dy)\nonumber\\
&=&\int_{[0,x-x_k]}\overline {F^{*(k-1)}}(x-y)F(dy)+\overline {F^{*(k-1)}}(x_k)\overline {F}(x-x_k)+
\int_{[0,x_k]}\overline {F}(x-y)F^{*(k-1)}(dy)\nonumber\\
&\leq&\varepsilon\int_{[0,x-x_k]}\overline {F^{*k}}(x-y)F(dy)+\overline {F^{*(k-1)}}(x_k)\overline {F}(x-x_k)
+\varepsilon\int_{[0,x_k]}\overline {F^{*2}}(x-y)F^{*(k-1)}(dy)\nonumber\\
&\leq&\varepsilon\int_{[0,x-x_k]}\overline {F^{*k}}(x-y)F(dy)+\varepsilon\overline {F^{*2}}(x)
+\varepsilon\int_{[x-x_k,x]}\overline {F^{*(k-1)}}(x-y)F^{*2}(dy)\nonumber\\
&\leq&2\varepsilon\overline {F^{*(k+1)}}(x).
\end{eqnarray*}
By the arbitrariness of $\epsilon$, we immediately obtain (\ref{191}).

In the following, we continue to prove Corollary \ref{Corol201} in Case (2). By (\ref{191}), we have $\overline{F}(x)=o(\overline{F^{*n}}(x))$. Furthermore, by $F^{*n}\in\mathcal{L}(\gamma)$ and Theorem 3 in \cite{EG1980}, $F^{*(n+1)}\in\mathcal{L}(\gamma)$.\hfill$\Box$
\vspace{0.2cm}

Now, we give three results on the closedness under the random convolution or composite distribution. In the following, let $\tau$ be a random variable with distribution $G$ satisfying $G(\{k\})=P(\tau=k)=p_k$ for all nonnegative integer $k\ge0$, where $\sum_{k=0}^\infty p_k=1$. Then the distribution $F^\tau=\sum_{k=0}^\infty p_k F^{*k}$ is called a random convolution generated by the random variable $\tau$ and some distribution $F$.
\begin{thm}\label{thm202}
For any integer $n\ge2$ and $\gamma\ge0$, let $F$ be a distribution such that $F^{*n} \in \mathcal{L}(\gamma)$. Assume that condition (\ref{00}), or condition (\ref{1040}), is satisfied. Further, suppose that $p_k=0,\ k\ge n$ and $p_n>0$. Then $H_n=\sum_{k=1}^n p_iF^{*k}\in\mathcal{L}(\gamma)$.
\end{thm}

\proof According to $F^{*n}\in\mathcal{L}(\gamma)$, (\ref{00}) or (\ref{1040}), (\ref{000}) in Lemma \ref{lemma202} holds. Further, by (\ref{000}) and $p_n>0$, we know that, for any $t>0$,
\begin{eqnarray*}
&&|\overline {H_n}(x-t)-e^{\gamma t}\overline {H_n}(x)|/\overline {H_n}(x)
\le\sum_{k=1}^n p_k|\overline{F^{*k}}(x-t)-e^{\gamma t}\overline{F^{*k}}(x)|/\big(\overline {F^{*n}}(x)p_n\big)\to0,
\end{eqnarray*}
thus $H_n\in\mathcal{L}(\gamma)$.\hfill$\Box$
\vspace{0.2cm}

\begin{thm}\label{thm203}
Let $F$ be a distribution such that $F^{*n}\in\cal{L}(\gamma)$ for some integer $n\ge1$ and some $\gamma\ge0$. Assume that $P(\tau\ge n)>0$ and the following condition is satisfied: for any $0<\varepsilon<1$, there is an integer $M=M(F,\varepsilon)$ large enough such that
\begin{eqnarray}\label{21}
\sum_{k=M}^\infty p_{k+1}\overline{F^{*k}}(x)\le\varepsilon\overline{F^{*\tau}}(x),\ for\ all\ x\ge0.
\end{eqnarray}
Further, suppose that the condition (\ref{00}) or condition (\ref{1040}) is satisfied. Then the random convolution $F^\tau\in\cal{L}(\gamma)$ {and
\begin{eqnarray}\label{210}
\liminf\overline{F^{*\tau}}(x)/\overline{F^{*n}}(x)\ge\sum_{m=1}^\infty m\sum_{k=nm}^{n(m+1)-1}p_k m(F^{*n})^{m-1}.
\end{eqnarray}
In particular, if $m(F^{*n})=\infty$, then $\overline{F^{*n}}(x)=o(\overline{F^{*\tau}}(x))$.}
\end{thm}
\proof According to (\ref{00}) and (\ref{21}), combined with Theorem \ref{thm202} {of the present paper} and the proof of Proposition 6.1 of Watanabe and Yamamuro (2010) \cite{WY2010}, we can get $F^\tau\in\cal{L}(\gamma)$. {Further, by Fatou's lemma, (\ref{210}) holds.} \hfill$\Box$
{\begin{thm}\label{thm204}
Let $F$ be a distribution such that $F^{*n}\in\cal{L}(\gamma)\cap\mathcal{OS}$ for some integer $n\ge1$ and some $\gamma\ge0$. Let $\tau$ be a random variable as in Theorem \ref{thm203} such that for some $\varepsilon_0>0$,
\begin{eqnarray}\label{23}
\sum_{m=1}^\infty\Big(\sum_{k=(m-1)n+1}^{mn}p_k\Big)\big(C^*(F^{*n})-m(F^{*n})+\varepsilon_0\big)^m<\infty.
\end{eqnarray}
Further, suppose that condition (\ref{00}) or condition (\ref{1040}) is satisfied. Then $F^{*\tau}\in\cal{L}(\gamma)\cap\mathcal{OS}$ and
\begin{eqnarray}\label{24}
\limsup\overline{F^{*\tau}}(x)/\overline{F^{*n}}(x)\le
\sum_{m=1}^\infty\Big(\sum_{k=(m-1)n+1}^{mn}p_k\Big)\sum_{i=0}^{m-1}\big(m(F^{*n})\big)^i\big(C^*(F^{*n})-m(F^{*n})\big)^{m-1-i}.
\end{eqnarray}
\end{thm}}
{The proof of the theorem is similar to that of Theorem 2 (2b) in \cite{XFW2015}, we omit the details of it.}

\begin{remark}\label{rem205}
Theorem \ref{thm203} {and Theorem \ref{thm204} in the case of $n=1$ are due to Lemma 4 and Corollary 1} of Yu and Wang \cite{YW2014}, which slightly improve Proposition 6.1 in \cite{WY2010}. However, {in this paper,} we do not require $F\in\cal{L}(\gamma)$.
In fact, if $F\in\mathcal{F}_i(\gamma),\ i=1,2,3,4$, see Definition \ref{defin301} and \ref{defin302} below, then $F$ satisfies condition (\ref{00}) (namely (\ref{104})) or condition (\ref{1040}), and $F^{*2}\in\cal{L}(\gamma)$, but $F\notin\cal{L}(\gamma)$. Therefore,
Theorem \ref{thm203} is a extension of Lemma 6 in \cite{YW2014} and Proposition 6.1 in \cite{WY2010}, {in which the conclusion for $\gamma=0$ is due to Theorem 2.2 (2a) in \cite{XFW2015}}.
\end{remark}

\begin{remark}\label{rem206}
Specifically, if
\begin{eqnarray}\label{222}
p_k>0\ for\ all\ k\ge1\ and\ p_{k+1}/p_k\to0\ as\ k\to\infty,
\end{eqnarray}
for example, $p_k=e^{-\lambda}\lambda^k/k!,\ k=0,1,\cdot\cdot\cdot$, where $\lambda$
is a positive constant, then condition (\ref{21}) is satisfied for any distribution $F$.

In addition, for distribution $F$, if Kesten's inequality holds, namely there are two positive constants $C$ and $\alpha$ such that
\begin{eqnarray}\label{22}
\overline{F^{*k}(x)}\le Ce^{\alpha k}\overline{F}(x)\ \ for\ all\ k\ge1\ and\ x\ge0,
\end{eqnarray}
and if $\sum_{k=1}^\infty p_k e^{\alpha k}<\infty,$ for example,
$p_k=qp^k,\ k=0,1,\cdot\cdot\cdot$, where $p$ and $q$ are two positive constants such that $p+q=1$ and $pe^\alpha<1$, then condition (\ref{21}) is satisfied. And, using the method applied in the proof of Lemma 5 in \cite{YW2014}, we can know (\ref{22}) holds for $F$ satisfying $F^{*n}\in\mathcal{L}(\gamma)\cap\mathcal{OS}$ for some $n\ge2$ and $\gamma\ge0$.
\end{remark}

\section{On the closedness under the convolution roots}
\setcounter{thm}{0}\setcounter{Corol}{0}\setcounter{lemma}{0}\setcounter{pron}{0}\setcounter{equation}{0}\setcounter{remark}{0}\setcounter{exam}{0}\setcounter{property}{0}\setcounter{defin}{0}

On the other hand, we also want to know that, under what conditions, the class $\mathcal{L}(\gamma)$ with some $\gamma>0$, or the class $\mathcal{L}(\gamma)\cap\mathcal{OS}$ is closed under convolution roots. Here, we give a initial result with its corollary, which also plays a role in the proof of Theorem \ref{thm101}.
\begin{pron}\label{pron101}
(1) Let $F$ be a distribution such that $F\in\mathcal{OS}$ and $F^{*2}\in\mathcal{L}(\gamma)$ for some $\gamma\ge0$. Further, assume that the condition (\ref{104}) or (\ref{1040}) is satisfied,
%\begin{eqnarray}\label{104} \liminf\overline F(x-t)/\overline F(x)\ge e^{\gamma t}. \end{eqnarray}
%or, in the case $\gamma>0$, the condition
%\begin{eqnarray}\label{1040}
%\overline{F}(x)=o(\overline{F^{*2}}(x))
%\end{eqnarray}
then $F\in\mathcal{L}(\gamma)$.

(2) Let $F$ be a distribution satisfying the condition (\ref{104}) for some $\gamma\ge0$. Further, assume that $F^{*2}\in\mathcal{L}(\gamma)\cap\mathcal{OS}$ and that, only when $\gamma=0$,
\begin{eqnarray}\label{105}
C^*(F^{*2})-2m(F^{*2})<1.
\end{eqnarray}
Then $F\in\mathcal{L}(\gamma)\cap\mathcal{OS}$.
\end{pron}
%(3) In (2), when $\gamma>0$, the condition (\ref{105}) can be canceled. That is, under the condition (\ref{104}) and $\gamma>0$, $F\in\mathcal{L}(\gamma)\cap\mathcal{OS}$ if and only if $F^{*2}\in\mathcal{L}(\gamma)\cap\mathcal{OS}$.
%In Theorem \ref{thm101}, the result in the case $\gamma>0$ is similar to the result in the case $\gamma=0$. However, the following result indicates that there is a substantial difference between heavy-tailed distribution and light-tailed distribution. We recall Theorem 1.2 (1) of \cite{XFW2015} that, there is a distribution $F$, which is neither long-tailed nor generalised subexponential, while $F^{*k}\in \mathcal{L}\cap\mathcal{OS}$ for all $k\ge 2$. But, this phenomenon can not happen in the case $\gamma>0$. It should be said that this difference is also surprising.
\begin{remark}\label{rem102}
The condition (\ref{105}) is necessary in some certain sense, see Theorem 2.2 (1) in \cite{XFW2015}. There, $F\in\mathcal{OL}\setminus\big(\mathcal{L}\cup\mathcal{OS}\big)$, $F^{*2}\in\mathcal{L}\cap\mathcal{OS}$, while the condition  (\ref{105}) does not hold,   otherwise, by Proposition \ref{pron101}, $F\in\mathcal{L}\cap\mathcal{OS}.$
\end{remark}
\proof (1) Since $F \in \mathcal{OS}$, $C^*(F)<\infty$. So, by $F^{*2} \in \mathcal{L}(\gamma)$, condition (\ref{104}) and Lemma \ref{lemma202} (or condition (\ref{1040})), we have,
\begin{eqnarray*}
|\overline F(x-t)-e^{\gamma t}\overline F(x)|/\overline F(x)\le C^*(F)|\overline F(x-t)-e^{\gamma t}\overline F(x)|/\overline {F^{*2}}(x)\to0,
\end{eqnarray*}
that is $F \in \mathcal{L}(\gamma)$.

(2) We first give the following two lemmas.
\begin{lemma}\label{lemma501}
Let $F$ be a distribution supported on $[-c,\infty)$ for some $c\ge0$ such that $m(F)<\infty$ and $F^{*2}\in\mathcal{L}(\gamma)$ for some $\gamma\ge0$. Then $F\in\mathcal{OS}$ if and only if

\begin{eqnarray}\label{501}
a(F)=\limsup_{k\to\infty}\limsup\int_{(k,x-k]}\overline{F}(x-y)F(dy)/\overline{F^{*2}}(x)=\limsup_{k\to\infty}a_k(F)<1.
\end{eqnarray}
\end{lemma}

\proof If (\ref{501}) holds, firstly, we are aim to prove $F\in\mathcal{OS}$
for the case $c=0$. For any $0<\varepsilon<\big(1-a(F)\big)/4$, there is an integer $k_0\ge1$ such that $a_k(F)<a(F)+\varepsilon$ and $\overline{F}(k)e^{\gamma k}<\varepsilon$ for all $k\ge k_0$. Thus, by $F^{*2}\in\mathcal{L}(\gamma)$, we have

\begin{eqnarray}\label{502}
\liminf\frac{2\int_{[0,k]}\overline{F}(x-y)F(dy)}{\overline{F^{*2}}(x)}
&=&\liminf\Big(1-\frac{\int_{(k,x-k]}\overline{F}(x-y)F(dy)}{\overline{F^{*2}}(x)}-\frac{\overline{F}(x-k)\overline{F}(k)}{\overline{F^{*2}}(x)}\Big)\nonumber\\
&\ge&1-\big(a(F)+\varepsilon\big)-\overline{F}(k)e^{\gamma k}\nonumber\\&>&\big(1-a(F)\big)/2>0.
\end{eqnarray}
Still by $F^{*2}\in\mathcal{L}(\gamma)$, the following inequality holds:
\begin{eqnarray}\label{503}
\liminf2\int_{[0,k]}\overline{F}(x-y)F(dy)/\overline{F^{*2}}(x)&\le&\liminf2\overline{F}(x-k)F(k)/\overline{F^{*2}}(x)\nonumber\\
&\le&2F(k)e^{\gamma k}/C^*(F).
\end{eqnarray}
According to (\ref{502}) and (\ref{503}), we know $C^*(F)<\infty$, thus $F\in\mathcal{OS}$.

Now we deal with the case $c>0$. Let $X$ be a random variable with distribution $F$ supported on $[-c,\infty)$. And let $Y=X+c$ be another random variable with a distribution $G$ supported on $[0,\infty)$. Clearly, $m(G)=e^{c\gamma}m(F)<\infty$ and $\overline{G^{*2}}(x)=\overline{F^{*2}}(x-2c)\sim e^{2c\gamma}\overline{F^{*2}}(x)$. Thus $G^{*2}\in\mathcal{L}(\gamma)$. Further, from
\begin{eqnarray*}
a(G)&=&\limsup_{k\to\infty}\limsup\int_{(k,x-k]}\overline{G}(x-y)G(dy)/\overline{G^{*2}}(x)\nonumber\\
&=&\limsup_{k\to\infty}\limsup\int_{(k-c,x-2c-(k-c)]}\overline{F}(x-2c-y)F(dy)/\overline{F^{*2}}(x-2c)=a(F)<1
\end{eqnarray*}
and the conclusion in the case $c=0$, we know $G\in\mathcal{OS}$. Finally, by the following fact
$$\overline{F}(x)\le\overline{G}(x)\le\overline{G^{*2}}(x)\sim e^{c\gamma}\overline{G^{*2}}(x+c)=O\big(\overline{G}(x+c)\big)=O\big(\overline{F}(x)\big),$$
we have $F\in\mathcal{OS}$.

If $F\in\mathcal{OS}$, we are aim to prove $a(F)<1$. Otherwise, we set $a(F)=1$, then from the following inequality
\begin{eqnarray*}
1&\ge&\liminf_{k\to\infty}\liminf\int_{[0,k]}\overline{F}(x-y)F(dy)/\overline{F^{*2}}(x)
+\limsup_{k\to\infty}\limsup\int_{(k,x-k]}\overline{F}(x-y)F(dy)/\overline{F^{*2}}(x)\nonumber\\
&\ge&\liminf_{k\to\infty}\liminf\overline{F}(x)F(k)/\overline{F^{*2}}(x)+1\nonumber\\
&=&1/C^*(F)+1,
\end{eqnarray*}
we know $C^*(F)=\infty$, which is contradictory to the fact that $F\in\mathcal{OS}$.\hfill$\Box$

\begin{lemma}\label{lemma502}
Let $F$ be a distribution such that $F^{*2}\in\mathcal{L}(\gamma)\cap\mathcal{OS}$ for some $\gamma>0$.
Further, assume that the conditions (\ref{104}) and (\ref{105}) are satisfied. Then $F\in\mathcal{OS}$.
\end{lemma}

\proof First, we prove $F\in\mathcal{OS}$ for the case that $\gamma\ge0$ and
$C^*(F^{*2})-2m(F^{*2})<1$.
By the partial integration method, we have
\begin{eqnarray*}
&&a(F)\le\limsup_{k\to\infty}\limsup\int_{(k,x-k]}\overline{F^{*2}}(x-y)F(dy)/\overline{F^{*2}}(x)\nonumber\\
&\le&\limsup_{k\to\infty}\limsup\Big(\overline{F^{*2}}(x-k)\overline{F}(k)+\int_{(k,x-k]}\overline{F^{*2}}(x-y)F^{*2}(dy)\Big)
/\overline{F^{*2}}(x)\nonumber\\
&=&C^*(F^{*2})-2m(F^{*2})<1.
\end{eqnarray*}
Thus, by Lemma \ref{lemma501}, $F\in\mathcal{OS}$.

Now, we deal with the case that $\gamma>0$ and $C^*(F^{*2})-2m(F^{*2})\ge1$. Let $X$ be a random variable with a distribution $F$. Then for any fixed $c>0$, the random variable $Y=X-c$ is a distribution $G$ supported on $[-c,\infty)$ such that $G^{*2}\in\mathcal{L}(\gamma)\cap\mathcal{OS}$, $\overline{G^{*2}}(x)\sim e^{-2c\gamma}\overline{F^{*2}}(x)$ and $m(G)=e^{-c\gamma}m(F)$.

According to Lemma 7 of Foss and Kirshunov \cite{FK2007} and the condition (\ref{104}), we know
\begin{eqnarray}\label{505}
\liminf\overline{F^{*2}}(x)/\overline{F}(x)\ge 2m(F).
\end{eqnarray}
And by Lemma 2 of Yu and Wang \cite{YW2014} and the properties mentioned above, we have
\begin{eqnarray}\label{506}
C^*(G^{*2})-2m(G^{*2})=e^{-2c\gamma}\big(C^*(F^{*2})-2m(F^{*2})\big).
\end{eqnarray}
By (\ref{505}), there are a constant $k_0>0$ such that for all $k\ge k_0$,
\begin{eqnarray*}
&&a_k(G)=\limsup\int_{(k,x-k]}\overline{G}(x-y)G(dy)/\overline{G^{*2}}(x)\nonumber\\
&\le&\limsup\int_{(k,x-k]}\overline{G^{*2}}(x-y)G(dy)/\big(4m(G)\overline{G^{*2}}(x)\big).
\end{eqnarray*}
Further, by (\ref{506}), there are a constant $c_0>0$ such that for all $c\ge c_0$,
\begin{eqnarray*}
a(G)&\le&\limsup_{k\to\infty}\limsup\big(\overline{G^{*2}}(x-k)\overline{G}(k)+\int_{(k,x-k]}\overline{G^{*2}}(x-y)G^{*2}(dy)\big)
/\big(16m^2(G)\overline{G^{*2}}(x)\big)\nonumber\\
&=&\big(C^*(G^{*2})-2m(G^{*2})\big)/\big(16e^{-2c\gamma}m^2(F)\big)\nonumber\\
&=&\big(C^*(F^{*2})-2m(F^{*2})\big)/\big(16m^2(F)\big)<1.
\end{eqnarray*}
Therefore, by Lemma \ref{lemma501}, $G\in\mathcal{OS}$, that is $F\in\mathcal{OS}$. \hfill$\Box$

Therefore, by $F\in\mathcal{OS},\ F^{*2}\in\mathcal{L}(\gamma)$, condition (\ref{104}) and (1) of the theorem, $F\in\mathcal{L}(\gamma)$.\hfill$\Box$

\begin{Corol}\label{Corol102}
Let $F$ be a distribution such that $F\in\mathcal{OS}\backslash\mathcal{L}(\gamma)$ for some $\gamma\ge0$. Further, assume that the condition (\ref{104}) is satisfied, then $F^{*k}\notin\mathcal{L}(\gamma)$ for all $k\ge2$.
\end{Corol}
\proof We assume that there exists some integer $n\ge2$ such that $F^{*n}\in\mathcal{L}(\gamma)$ and $F^{*i}\notin\mathcal{L}(\gamma)$ for all $1\le i<n$.

When $\gamma>0$, under condition (\ref{104}), by $F^{*n}\in\mathcal{L}(\gamma)$ and Corollary \ref{Corol201}, we can get $F^{*i}\in\mathcal{L}(\gamma)$ for all $i\ge n$. Hence $F^{*2(n-1)}\in\mathcal{L}(\gamma)$, but $F^{*(n-1)}\notin\mathcal{L}(\gamma)$. From (1) of the theorem, we know $F^{*(n-1)}\notin\mathcal{OS}$, which is  contradictory to $F\in\mathcal{OS}$.

When $\gamma=0$, the conclusion comes from Remark \ref{rem204}.\hfill$\Box$

\section{A transformation}
\setcounter{thm}{0}\setcounter{Corol}{0}\setcounter{lemma}{0}\setcounter{pron}{0}\setcounter{equation}{0}
\setcounter{remark}{0}\setcounter{exam}{0}\setcounter{property}{0}\setcounter{defin}{0}

In order to prove Theorem 1.1, we need to construct some light-tailed distributions with some interesting properties, some of which come from the corresponding heavy-tailed distributions through a transformation between the distributions. For some constant $\gamma>0$ and distribution $F_0$, we define the distribution $F_\gamma$ in the form
\begin{eqnarray}\label{303}
\overline{F_\gamma}(x)=\textbf{\emph{\emph{1}}}(x<0)+e^{-\gamma x}\overline{F_0}(x)\textbf{\emph{\emph{1}}}(x\ge0),\ x\in(-\infty,\infty).
\end{eqnarray}
Clearly, $F_\gamma$ is light-tailed. In this way, {we can characterize corresponding light-tailed distribution $F_\gamma$ through certain heavy-tailed distribution $F_0$ with some good properties}, see Kl\"{u}ppelberg \cite{K1989}, Xu et al. \cite{XSW2014}, and so on.
Here, we give the following new properties of the distribution $F_\gamma$.
\begin{lemma}\label{lemma301}
For some $\gamma>0$, let $F_0,F_\gamma;F_{i0},F_{i\gamma}, i = 1,2$ be three pairs of distributions defined in (\ref{303}), respectively. Then the following conclusions hold.\\
(1) For all $t>0$ and $x\ge t$,
\begin{eqnarray*}\label{3031}
\overline{F_\gamma}(x-t)/\overline{F_\gamma}(x)\ge e^{\gamma t}.
\end{eqnarray*}
(2) For $i=1,2$, if
\begin{eqnarray}\label{3032}
\mid\overline{F_{i0}}(x-t)-\overline{F_{i0}}(x)\mid=o\big(\overline{F_{10}*F_{20}}(x)\big)\ for\ all\ t>0,
\end{eqnarray}
then
\begin{eqnarray*}\label{3033}
\mid\overline{F_{i\gamma}}(x-t)-e^{\gamma t}\overline{F_{i\gamma}}(x)\mid=o\big(\overline{F_{1\gamma}*F_{2\gamma}}(x)\big)\ for\ all\ t>0.
\end{eqnarray*}
(3) For $1\le i\neq j\le 2$, if ${\mu(F_{j0})}=\int_0^\infty \overline{F_{j0}}(y)dy=\infty$, then
\begin{eqnarray}\label{3034}
\overline{F_{i\gamma}}(x)=o\big(\overline {F_{i\gamma}*F_{j\gamma}}(x)\big).
\end{eqnarray}
Particularly, if $F_{10}=F_{20}=F_0\in\cal{OS}$,
%then $\overline{F_{\gamma}}(x)=o\big(\overline {F_{\gamma}^{*2}}(x)\big)$.
 then
\begin{eqnarray}\label{30340}
\overline{F_{\gamma}^{*2}}(x)\sim \gamma e^{-\gamma x}\int_0^x\overline{F_{0}}(x-y)\overline{F_{0}}(y)dy
=\gamma \int_0^x\overline{F_{\gamma}}(x-y)\overline{F_{\gamma}}(y)dy.
\end{eqnarray}
\end{lemma}

\proof (1) This is an obvious fact.\\
(2) For $1\le i\neq j\le 2$, according to (\ref{303}),
\begin{eqnarray}\label{3036}
\overline {F_{i\gamma}*F_{j\gamma}}(x)&=&\overline {F_{i\gamma}}(x)+\int_{[0,x]}\overline {F_{j\gamma}}(x-y)F_{i\gamma}(dy)\nonumber\\
&=&e^{-\gamma x}\Big(\overline {F_{i0}*F_{j0}}(x)
+\gamma\int_{[0,x]}\overline {F_{j0}}(x-y)\overline F_{i0}(y)dy\Big).
\end{eqnarray}
Further, by (\ref{3032}), we have
\begin{eqnarray*}
\mid\overline F_{i\gamma}(x-t)-e^{\gamma t}\overline F_{i\gamma}(x)\mid/\overline {F_{i\gamma}*F_{j\gamma}}(x)
\leq e^{\gamma t}\mid\overline F_{i0}(x-t)-\overline F_{i0}(x)\mid/\overline {F_{i0}*F_{j0}}(x)\to0,
\end{eqnarray*}
for $1\le i\neq j\le 2$.\\
(3) According to (\ref{3036}), we know that
$$\overline {F_{i\gamma}*F_{j\gamma}}(x)\ge e^{-\gamma x}\gamma\overline {F_{i0}}(x)\int_{[0,x]}\overline {F_{j0}}(x-y)dy
=\gamma\overline {F_{i\gamma}}(x)\int_{[0,x]}\overline {F_{j0}}(y)dy.$$
Therefore, by $F_0\in\cal{OS}$, (\ref{3036}) and $\mu(F_{0})=\infty$, (\ref{3034}) holds. Further, (\ref{30340}) is proved.\hfill$\Box$

\begin{lemma}\label{lemma302}
For $i=1,2$, let $F_{i0}$ be an absolutely continuous distribution with density $f_{i0}$ such that
\begin{eqnarray}\label{l30101}
F_{i0}(x-t,x]=\overline{F_{i0}}(x-t)-\overline{F_{i0}}(x)=O(f_{i0}(x-t)+f_{i0}(x))\ for\ all\ 0\le t< x
\end{eqnarray}
and
\begin{eqnarray}\label{l3011}
\int_{[x/2,x]}\overline F_{i0}(x-y)F_{j0}(dy)=o\Big(\int_{[x/2,x]}\overline F_{i0}(x-y)\overline F_{j0}(y)dy\Big),\ 1\le j\not=i\le2,
\end{eqnarray}
then for any $\gamma>0,\ F_{1\gamma}*F_{2\gamma}\in\mathcal{L}(\gamma)$.
\end{lemma}

\begin{remark}\label{rem301}
Firstly, in the lemma, $F_{1\gamma}$ or $F_{2\gamma}$ is not required to belong to the class $\mathcal{L}(\gamma)$, which is a difference compared to Theorem 3 {in} \cite{EG1980}, Theorem 1.1 (1b) {in} \cite{XFW2015} and Theorem \ref{thm201} of the paper.

Secondly, if $F_{10}=F_{20}=F_0\in\mathcal{L}$, then it is clear that $F_{1\gamma}=F_{2\gamma}=F_\gamma$ and $F_\gamma^{*2}\in\mathcal{L}(\gamma)$ for any $\gamma>0$.
%In this case, the above condition (\ref{l3010}) or condition (\ref{l301}) is not necessary, for example, in Example 2.2 of \cite{WXCY2014}, $f_0(x)=o(\overline{F_0}(x))$, then (\ref{l301}) is satisfied, while (\ref{l3010}) is not holds.
More interestingly, there are some distributions $F_0\notin\mathcal{L}$ satisfying (\ref{l30101}) and (\ref{l3011}), for example, see distributions in the classes $\mathcal{F}_i(0)$ for $i=1,2,3$, which have appeared in Section 5 below. Thus $F_\gamma\notin\mathcal{L}(\gamma)$, while $F_\gamma^{*2}\in\mathcal{L}(\gamma)$ by Lemma \ref{lemma302}.
%So we say that, as Theorem 3 of Embrechts and Goldie \cite{EG1980} and Theorem 1.1 (1b) of Xu et al. \cite{XFW2015}, the lemma is an important result to prove $F_\gamma^{*2}\in\mathcal{L}(\gamma)$.

%On the contrary, if $F_0\in\mathcal{OS}^*\setminus\mathcal{S}$, then (\ref{l301}) does not hold, even $F_0\in\mathcal{L}$, see Example 2.2 in \cite{WXCY2014}.

Finally, we note that conditions (\ref{l30101}) and (\ref{l3011}) can not be deduced from each other, {and that} there is a distribution $F$ belonging to the class $\mathcal{L}(\gamma)$ but not satisfying condition (\ref{l30101}), see Subsection 6.3 for the details.
\end{remark}

\proof We note that $F_{10}$, $F_{20}$ are absolutely continuous. Thus, from (\ref{3036}), we have
%\begin{eqnarray}\label{exam203-21}
$$\overline {F_{1\gamma}*F_{2\gamma}}(x)=e^{-\gamma x}S(x)$$
%\end{eqnarray}
 for all $x\ge0$, where
\begin{eqnarray*}
S(x)=\overline {F_{10}*F_{20}}(x)+2\gamma \int_{[0,x/2]}\overline{F_{20}}(x-y)\overline{F_{10}}(y)dy=\overline {F_{10}*F_{20}}(x)+2\gamma T(x).
\end{eqnarray*}

Hence $F_{1\gamma}*F_{2\gamma}\in\mathcal{L}(\gamma)$ is equivalent to that the function $S$ belongs to the long-tailed function class
$$\mathcal{L}_d=\{f:f(x)>0\ for\ all\ x\in\ (-\infty,\infty)\ and\ f(x-t)\sim f(x)\ for\ any\ t\in(-\infty,\infty)\}.$$

For any $t>0$, on the one hand, we have
\begin{eqnarray*}
&&S(x-t)-S(x)=\overline {F_{10}*F_{20}}(x-t)-\overline {F_{10}*F_{20}}(x)+\gamma \int_{[0,(x-t)/2]}\overline{F_{20}}(x-t-y)\overline{F_{10}}(y)dy\nonumber\\
&&\ \ \ \ -\gamma\int_{[0,x/2]}\overline{F_{20}}(x-y)\overline{F_{10}}(y)dy
+\gamma\int_{[0,(x-t)/2]}\overline{F_{10}}(x-t-y)\overline{F_{20}}(y)dy
-\gamma\int_{[0,x/2]}\overline{F_{10}}(x-y)\overline{F_{20}}(y)dy\nonumber\\
&\ge&\gamma \int_{[0,(x-t)/2]}(\overline{F_{20}}(x-t-y)-\overline{F_{20}}(x-y))\overline{F_{10}}(y)dy-\gamma \int_{((x-t)/2,x/2]}\overline{F_{20}}(x-y)\overline{F_{10}}(y)dy\nonumber\\
& &\ \ \ \ +\gamma \int_{[0,(x-t)/2]}(\overline{F_{10}}(x-t-y)-\overline{F_{10}}(x-y))\overline{F_{20}}(y)dy-\gamma \int_{((x-t)/2,x/2]}\overline{F_{10}}(x-y)\overline{F_{20}}(y)dy\nonumber\\
&\geq&-\gamma\int_{[(x-t)/2,(x+t)/2]}\overline{F_{20}}(x-y)\overline{F_{10}}(y)dy\nonumber\\
&\geq&-\gamma t\overline{F_{10}}((x-t)/2)
\overline{F_{20}}((x-t)/2)\nonumber\\
&\geq&-\gamma t\Big(\int_{[(x-t)/2,x-t]}\overline F_{20}(x-t-y)F_{10}(dy)+\overline F_{20}((x-t)/2)\overline F_{10}(x-t)\Big).
\end{eqnarray*}
Thus, by (\ref{l3011}) we know that
\begin{eqnarray}\label{exam203-31}
\liminf\big(S(x-t)-S(x)\big)/S(x-t)\ge0.
\end{eqnarray}

On the other hand, by (\ref{l30101}), we have

\begin{eqnarray}\label{exam203-41}
&&S(x-t)-S(x)\leq\overline {F_{10}*F_{20}}(x-t)-\overline {F_{10}*F_{20}}(x)\nonumber\\
&&\ \ \ \ \ \ \ \ +\gamma \int_{[0,(x-t)/2]}(\overline{F_{20}}(x-t-y)-\overline{F_{20}}(x-y))\overline{F_{10}}(y)dy\nonumber\\
&&\ \ \ \ \ \ \ \ +\gamma \int_{[0,(x-t)/2]}(\overline{F_{10}}(x-t-y)-\overline{F_{10}}(x-y))\overline{F_{20}}(y)dy\nonumber\\
&=&\sum_{i=1}^{2}(\overline {F_{i0}}(x-t)-\overline {F_{i0}}(x))+\int_{[(x-t)/2,x-t]}\overline F_{20}(x-t-y)F_{10}(dy)+\int_{[(x-t)/2,x-t]}\overline F_{10}(x-t-y)F_{20}(dy)\nonumber\\
& &\ \ \ \ \ \ \
+O\Big(\sum_{1\le i\not=j\le 2}\int_{[0,(x-t)/2]}\big(f_{i0}(x-t-y)+f_{i0}(x-y)\big)\overline{F_{j0}}(y)dy\Big)\nonumber\\
&=&O\Big(\sum_{1\le i\not=j\le 2}\Big(\int_{[x-t,x]}\overline F_{i0}(x-y)F_{j0}(dy)+\int_{[(x-t)/2,x-t]}\overline F_{i0}(x-t-y)F_{j0}(dy)\nonumber\\
&&\ \ \ \ \ \ \ \ + \int_{[(x+t)/2,x]}\overline{F_{i0}}(x-y)F_{j0}(dy)\Big)\Big)\nonumber\\
&=&O\Big(\sum_{1\le i\not=j\le 2}\Big(\int_{[x/2,x]}\overline F_{i0}(x-y)F_{j0}(dy)+\int_{[(x-t)/2,x-t]}\overline F_{i0}(x-t-y)F_{j0}(dy)\Big)\Big)\nonumber\\
&=&o\big(S(x)+S(x-t)\big).
\end{eqnarray}
From (\ref{exam203-31}), we know that, there is a positive constant $x_0=x_0(S,t)$ such that for all $x\ge x_0$,
$$S(x)\le 2S(x-t).$$
Thus, by (\ref{exam203-41}), we have
\begin{eqnarray}\label{exam203-51}
{\limsup \big(S(x-t)-S(x)\big)/S(x-t)\le0.}
\end{eqnarray}

Combining with (\ref{exam203-31}) and (\ref{exam203-51}), we know that $T\in\mathcal{L}_d$.\hfill$\Box$

%Next we prove $F_{1\gamma}*F_{2\gamma}\in\mathcal{OS}$. It suffices to prove. By (\ref{l3011}) it is easy to verify $T(x)=O\Big(\int_{0}^{x}\overline F_{10}(x-y)\overline F_{20}(y)dy\Big)=:O(g(x))$.

\section{Proof of Theorem \ref{thm101}}
\setcounter{thm}{0}\setcounter{Corol}{0}\setcounter{lemma}{0}\setcounter{pron}{0}\setcounter{equation}{0}
\setcounter{remark}{0}\setcounter{exam}{0}\setcounter{property}{0}\setcounter{defin}{0}

In order to find more L$\acute{e}$vy spectral distributions satisfying the requirement of Theorem \ref{thm101}, we construct the following four distribution classes with different properties, see Definitions \ref{defin301} and \ref{defin302} below. To this end, we first recall two classes of heavy-tailed distributions introduced by \cite{XFW2015}, then construct a new heavy-tailed distribution class.

Let $\alpha\in[1/2,1)$, $r=1+1/\alpha$, $b\ge1$ and $s\ge1$ be constants.
Let $a>1$ be large enough such that $a^{r} > 2^{s+2}a$, and let sequence
$A=\{a_n\}$ be given by $a_n =
a^{r^{n}}$ for $n = 0,1,\cdots.$ Let $\eta$ be a
discrete random variable distributed by $\textbf{P}(\eta=a_n)=Ca_n^{-\alpha}$, where
$C=(\sum_{0}^{\infty}a_n^{-\alpha})^{-1}$ is the normalising constant.
Let $U$ be a random variable having uniform distribution in the
interval $(0,1)$. Further, suppose that $U$ and $\eta$ are independent to each other.
Let $\mathcal{F}_1(0)$ be the class consisting of  4-parametric distributions
$F_0=F_0(\alpha,a,b,s)$ of random variables
\be \label{301}
\xi=\eta(1+U^{1/b})^s,
\ee which are absolutely continuous with densities $f_0=f_0(\alpha,a,b,s)$.

Next, {let $s\in(1,2)$ and $\alpha\in(1-s^{-1},s^{-1})$ be
constants. Further, assume that the constants $a$ and $r$, the sequence $A=\{a_n\}$ and the
random variables $\eta$ and $U$ are defined as before.}
Let $\mathcal{F}_2(0)$ be the class consisting of 3-parametric  heavy-tailed distributions
$F_0=F_0(\alpha,a,s)$ of random variables
\be \label{302}
\xi=\eta^{s^{-1}}(1+U)^{s^{-1}},
\ee
which are absolutely continuous with density $f_0=f_0(\alpha,a,s)$.

Finally, let $\alpha\in(3/2,(\sqrt{5}+1)/2)$ and $r=1+1/\alpha$ be constants.
{Assume $a>1$ is large enough such  that $a^{r} > 8a$. Further, suppose that the sequence $A=\{a_n\}$ and the
random variables $\eta$ and $U$ are defined as before.}
Let $\mathcal{F}_3(0)$ be the class consisting of 2-parametric  distributions
$F_0=F_0(\alpha,a)$  of random variables
\be \label{304}
\xi=\eta(1+U),
\ee
which are absolutely continuous with density $f_0=f_0(\alpha,a)$.

Obviously, these distributions of random variables $\xi$ in (\ref{301}), (\ref{302}) and (\ref{304}) are natural {and intuitive}. For example, in (\ref{301}), random variables $U^{1/b}$ and $\eta$ may be interpreted as the interest of a financial institution and the capital of an investor, respectively, then the random variable $\xi=\eta(1+U^{1/b})^{s}$ is the total income of the investor at time $s$.

As pointed out in \cite{XFW2015}, if $F_0\in\mathcal{F}_1(0)$,
then $F_0^{*k}\in\mathcal{L}$ for all $k\ge2$, while $F_0\in\mathcal{OL}\setminus\mathcal{L}$ with infinite mean; and if $F_0\in\mathcal{F}_2(0)$, then $F_0^{*k}\in\mathcal{L}$ for all $k\ge2$,
while $F_0\notin\mathcal{OL}$ with infinite mean. In the following, however, we will find that if
$F_0\in\mathcal{F}_3(0)$, then its mean $\mu(F_0)<\infty$ and $F_0^{*k}\notin\mathcal{L}$ for all $k\ge1$. In spite of this, the corresponding distribution $F_\gamma$ in (\ref{303}) has very good properties.
%Nevertheless, there is a striking fact that, in (\ref{303}), corresponding $F_{i\gamma}^{*n}\in\mathcal{L}(\gamma)$,
%further $F_{i\gamma}^{*\tau}\in\mathcal{L}(\gamma)$ under certain condition for $\tau$,
%while $F_{i\gamma}\notin\mathcal{L}(\gamma)$, for all $\gamma>0,\ n\ge2$ and $i=1,2,3$, see below for the details.

Based on (\ref{303}) and the classes $\mathcal{F}_i(0),i=1,2,3$, we define the following three light-tailed distribution subclasses, respectively.
\begin{defin}\label{defin301}
For $i=1,2,3$, we say that the distribution $F_\gamma$ for some constant $\gamma>0$
belongs to the class $\mathcal{F}_i(\gamma)$, if the corresponding distribution $F_0$ in (\ref{303})
belongs to the class $\mathcal{F}_i(0)$.
\end{defin}

For notational convenience, we replace $F_\gamma$ and $f_\gamma$ with $F$ and $f$ in the following text.
And in order to prove Theorem \ref{thm101}, we give some properties
of the classes $\mathcal{F}_i(\gamma),\ i=1,2,3$, respectively.

\begin{pron}\label{pron301}
If $F\in \mathcal{F}_1(\gamma)$ for some $\gamma>0$, then {conditions (\ref{104}) and (\ref{1040}) are satisfied, $m(F)=\infty$ and $F^{*k}\in \mathcal{L}(\gamma)\setminus\mathcal{OS}$ for all $k\ge 2$, while $F\in\cal{OL}\setminus\big(\mathcal{L}(\gamma)\cup\mathcal{OS}\big)$}. Further, if the condition (\ref{21}) is satisfied with some non-negative integer-valued
random variable $\tau$, then $F^{*\tau}\in \mathcal{L}(\gamma)\setminus\mathcal{OS}$.
\end{pron}

\proof According to $F_0\in\mathcal{OL}\setminus\mathcal{L}$, we know that $F\in\cal{OL}\setminus\mathcal{L}(\gamma)$.
From $\mu(F_0)=\infty$, we know that {$m(F)=\infty$, thus $F^{*k},k\ge1$ and $F^{*\tau}\notin \mathcal{OS}$.  And {conditions (\ref{104}) and (\ref{1040}) follow from Lemma \ref{lemma301} (1) and (3), respectively}. Therefore, by Corollary \ref{Corol201} (2) and Theorem \ref{thm203},
we only need to show that $F^{*2}$ belongs to $\mathcal{L}(\gamma)$.

From \cite{XSW2014}, we know that if $F_0\in\mathcal{F}_1(0)$, then
\begin{eqnarray}\label{exam301-6}
f_0(x)=Cbs^{-1}\sum_{n=0}^{\infty}x^{1/s-1}a_n^{-\alpha-1/s}
\big((xa_n^{-1})^{1/s}-1\big)^{b-1}\textbf{1}(x\in [a_n,2^sa_n))
\end{eqnarray}
and
\begin{eqnarray}\label{exam301-7}
\overline{F_0}(x)&=&\textbf{1}(x<
a_0)+\sum\limits_{n=0}^{\infty}\Big(\big(\sum\limits_{i=n}^{\infty}
Ca_i^{-\alpha}-Ca_n^{-\alpha}((x/a_n)^{1/s}-1)^{b}\big)
\textbf{1}(x\in[a_n,2^sa_n))\nonumber\\
& &
+\sum\limits_{i=n+1}^{\infty}Ca_i^{-\alpha}\textbf{1}(x\in[2^sa_n,
a_{n+1}))\Big).
\end{eqnarray}

For all $t>0$ and $n = 0,1,\cdot \cdot \cdot$, when $x\in[a_n,2^sa_n+t)$,
$$\overline F_0(x-t)-\overline F_0(x)\leq Cbts^{-1}a_n^{-\alpha-1}=O(f_0(x-t)+f_0(x));$$
when $x\in[2^sa_n+t,a_{n+1})$,
$$\overline F_0(x-t)-\overline F_0(x)=f_0(x-t)+f_0(x)=0,$$
thus (\ref{l30101}) holds. Hence by Lemma \ref{lemma302}, in order to prove $F^{*2}\in\mathcal{L}(\gamma)$, it suffices to prove (\ref{l3011}), that is
\begin{eqnarray}\label{exam301-70}
W(x)=\int_{[x/2,x]}\overline F_0(x-y)F_0(dy)=o\Big(\int_{[x/2,x]}\overline F_0(x-y)\overline F_0(y)dy\Big)=o(T(x)).
\end{eqnarray}

For all $n = 0,1,\cdot \cdot \cdot$, because $W(x)=0$ for $x\in[2^{s+1}a_n,a_{n+1})$, we only need to deal with $W(x)$ in the following two cases:
$i)\ x\in[a_n,2^{s+1}a_n-a_n^{5/6})$ and $ii)\ x\in[2^{s+1}a_n-a_n^{5/6},2^{s+1}a_n)$.

In case $i)$, by (\ref{exam301-6}), (\ref{exam301-7}) and $\overline{F_0}((x+2^{-1}a_n^{5/6})/2)\ge\overline{F_0}(2^{s}a_n-4^{-1}a_n^{5/6})$, we have,

\begin{eqnarray}\label{exam201-4}
W(x)/T(x)&\leq&\Big(Cbs^{-1}a_n^{-\alpha-1}\int_{[x/2,x]}\overline{F_0}(x-y)dy\Big)/
\Big(\int_{[x/2,(x+2^{-1}a_n^{5/6})/2]}\overline{F_0}(x-y)\overline{F_0}(y)dy\Big)\nonumber\\
&\leq&\Big(4Cbs^{-1}a_n^{-\alpha-1}\int_{[0,2^{s}a_n]}\overline{F_0}(y)dy\Big)/
\Big(a_n^{5/6}\overline{F_0}^2(2^{s}a_n-4^{-1}a_n^{5/6})\Big)\nonumber\\
&=&O(a_n^{-1/2})\to0,~~~n\rightarrow\infty.
\end{eqnarray}

Then, in case $ii)$, by (\ref{exam301-6}) and (\ref{exam301-7}), we have
$$\int_{[x/2,x]}\overline{F_0}(x-y)\overline{F_0}(y)dy\ge
Ca_n^{-\alpha-1}\int_{[0,x/2]}\overline{F_0}(y)dy.$$
Thus,
\begin{eqnarray}\label{exam201-5}
W(x)/T(x)&\leq&\int_{[x/2,2^{s}a_n]}\overline{F_0}(x-y)F_0(dy)/
\Big(\overline{F_0}(a_{n+1})\int_{[0,x/2]}\overline{F_0}(y)dy\Big)\nonumber\\
&\leq&\Big(bs^{-1}\int_{[2^{s}a_n-a_n^{5/6},2^{s}a_n]}\overline{F_0}(y)dy\Big)/
\int_{[a_n/2,a_n]}\overline{F_0}(y)dy\nonumber\\
&=&O(a_n^{-1/3})\to0,~~~n\rightarrow\infty.
\end{eqnarray}

According to (\ref{exam201-4}) and (\ref{exam201-5}), (\ref{exam301-70}) holds.
Therefore, $F^{*2}\in\cal{L}(\gamma).$\hfill$\Box$

\begin{pron}\label{pron302}
If $F\in \mathcal{F}_2(\gamma)$ for some $\gamma>0$, then {conditions (\ref{104}) and (\ref{1040}) are satisfied, $m(F)=\infty$ and $F^{*k}\in \mathcal{L}(\gamma)\setminus\mathcal{OS}$ for all $k\ge 2$, while $F\notin \mathcal{OL}$}. Further, if the condition (\ref{21}) is satisfied with some non-negative integer-valued random variable $\tau$, then $F^{*\tau}\in \mathcal{L}(\gamma )\setminus\mathcal{OS}$.
\end{pron}

\proof According to (\ref{303}) and $F_0\notin\cal{OL}$, we know $F\notin\cal{OL}$. Then by Lemma \ref{lemma301} and $\mu(F_0)=\infty$, we know that $F$ satisfies conditions (\ref{104}) and (\ref{1040}) and $m(F)=\infty$, thus $F^{*k},k\ge1$ and $F^{*\tau}\notin \mathcal{OS}$. So, for other conclusions, by Corollary \ref{Corol201} (2) and Theorem \ref{thm203}, we only need to show $F^{*2}\in\mathcal{L}(\gamma)$.

We have already known
\begin{eqnarray}\label{exam301-4}
f_0(x)=Cs \sum_{n=0}^{\infty}x^{s-1}a_n^{-\alpha-1}\textbf{1}\big(x\in
[a_n^{1/s},(2a_n)^{1/s})\big)
\end{eqnarray}
and
\begin{eqnarray}\label{exam301-5}
\overline{F_0}(x)&=&\textbf{1}(x< a_0)
+C\sum\limits_{n=0}^{\infty}\Big(\big(\sum\limits_{i=n}^{\infty}
a_i^{-\alpha}-a_n^{-\alpha}(a_n^{-1} x^s-1)\big)\textbf{1}(x\in
[a_n^{1/s},(2a_n)^{1/s}))\nonumber\\
&&+\sum\limits_{i=n+1}^{\infty}a_i^{-\alpha}\textbf{1}(x\in[(2a_n)^{1/s},a_{n+1}^{1/s}))\Big).
\end{eqnarray}

Similar to the proof of Proposition \ref{pron301}, we just need to prove (\ref{exam301-70}).

For all $n = 0,1,\cdot \cdot \cdot$, because $W(x)=\int_{[x/2,x]}\overline F_0(x-y)F_0(dy)=0$ for $x\in[2(2a_n)^{1/s},a_{n+1}^{1/s})$, we only need to deal with $W(x)$ in the following three cases:
$i)\ x\in[a_n^{1/s},2(2a_n)^{1/s}-a_n^{3 /4s})$; $ii)\ x\in[2(2a_n)^{1/s}-a_n^{3 /4s},2(2a_n)^{1/s}-a_n^{s^{-1}-4^{-1}})$ and $iii)\ x\in[2(2a_n)^{1/s}-a_n^{s^{-1}-4^{-1}},2(2a_n)^{1/s})$.

In case $i)$, by (\ref{exam301-4}), (\ref{exam301-5}) and $\overline{F_0}(x/2)\le\overline{F_0}((2a_n)^{1/s}-2^{-1}a_n^{3 /4s})$,  we have,

\begin{eqnarray}\label{exam202-4}
W(x)/T(x)&\leq&2Csa_n^{-\alpha-s^{-1}}\int_{[x/2,x]}\overline{F_0}(x-y)dy/
\Big(\overline{F_0}((2a_n)^{1/s}-2^{-1}a_n^{3 /4s})\int_{[x/2,(x+2^{-1}a_n^{3 /4s})/2]}\overline{F_0}(y)dy\Big)\nonumber\\
&\leq&8Csa_n^{-\alpha-7\cdot(4s)^{-1}}\int_{[0,(2a_n)^{1/s}-2^{-1}a_n^{3/4s}]}\overline{F_0}(y)dy/
\overline{F_0}^2((2a_n)^{1/s}-4^{-1}a_n^{3 /4s})\nonumber\\
&=&O(a_n^{-1/4s})\to0,~~~n\rightarrow\infty.
\end{eqnarray}

Then, in case $ii)$, from (\ref{exam301-4}), (\ref{exam301-5}) and $(2a_n)^{1/s}-2^{-1}a_n^{3/4s} \le x/2 \le y \le (2a_n)^{1/s}$,

\begin{eqnarray}\label{exam202-5}
&&W(x)/T(x)\leq\int_{[x/2,(2a_n)^{1/s}]}\overline{F_0}(x-y)F_0(dy)/
\int_{[x/2,(2a_n)^{1/s}]}\overline{F_0}(y)\overline{F_0}(x-y)dy\nonumber\\
&\leq&Csa_n^{-\alpha-(4s)^{-1}}\overline{F_0}(x-(2a_n)^{1/s}) /\Big(\overline{F_0}((2a_n)^{1/s}-2^{-1}a_n^{s^{-1}-4^{-1}})\nonumber\\
& &\cdot\int_{[(2a_n)^{1/s}-2^{-1}a_n^{s^{-1}-4^{-1}},
(2a_n)^{1/s}-4^{-1}a_n^{s^{-1}-4^{-1}}]}\overline{F_0}(y)dy\Big)\nonumber\\
&\leq&Csa_n^{-\alpha-5\cdot(4s)^{-1}+4^{-1}}\overline{F_0}((2a_n)^{1/s}-a_n^{3 /4s})/
\overline{F_0}^2((2a_n)^{1/s}-4^{-1}a_n^{s^{-1}-4^{-1}})\nonumber\\
&=&O(a_n^{-3\cdot4^{-1}(2s^{-1}-1)})\to0,~~~n\rightarrow\infty.
\end{eqnarray}

Finally, we consider case $iii)$. By (\ref{exam301-4}), (\ref{exam301-5}) and
$$\int_{[x/2,x]}\overline{F_0}(x-y)\overline{F_0}(y)dy\ge
Ca_n^{-\alpha-1}\int_{[0,x/2]}\overline{F_0}(y)dy,$$
we have
\begin{eqnarray}\label{exam202-6}
W(x)/T(x)&\leq&
\int_{[x/2,(2a_n)^{1/s}]}\overline{F_0}(x-y)F_0(dy)/
\Big(\overline{F_0}(a_{n+1}^{1/s})\int_{[0,x/2]}\overline{F_0}(y)dy\Big)\nonumber\\
&\leq&\Big(2sa_n^{1-s^{-1}}\int_{[(2a_n)^{1/s}-a_n^{s^{-1}-4^{-1}},(2a_n)^{1/s}]}\overline{F_0}(y)dy\Big)/
\int_{[a_n^{1/s}/2,a_n^{1/s}]}\overline{F_0}(y)dy\nonumber\\
&=&O(a_n^{2^{-1}-s^{-1}})\to0,~~~n\rightarrow\infty.
\end{eqnarray}

According to (\ref{exam202-4})-(\ref{exam202-6}), we get (\ref{exam301-70}).\hfill$\Box$
\vspace{0.2cm}

In the following proposition, we find a surprising phenomenon. There is a distribution $F_0$ with bad properties, but its corresponding distribution $F$ may still enjoy good ones.
\begin{pron}\label{pron303}
If $F\in \mathcal{F}_3(\gamma)$ for some $\gamma>0$, then {$m(F)<\infty$ and $F^{*k}\in \big(\mathcal{L}(\gamma )\cap\mathcal{OS}\big)\setminus\mathcal{S}(\gamma)$ for all $k\ge 2$, while $F\in\cal{OL}\setminus\big(\mathcal{L}(\gamma)\cup\mathcal{OS}\big)$ and $F_0^{*k}\in\cal{OS}\setminus\mathcal{L}$ for all $k\ge1$. In addition, for the distribution $F$,  condition (\ref{104}) is satisfied, but condition (\ref{1040}) is not}. Further, if the condition (\ref{21}) is satisfied with some non-negative integer-valued random variable $\tau$, then $F^{*\tau}\in {\big(\mathcal{L}(\gamma)\cap\mathcal{OS}\big)\setminus\mathcal{S}(\gamma)}$.
\end{pron}

\proof For any $F_0\in\mathcal{F}_3(0)$, it is easy to verify
\begin{eqnarray}\label{exam203-7}
f_0(x)=\sum_{n=0}^{\infty}Ca_n^{-\alpha-1}\textbf{1}(x\in [a_n,2a_n))
\end{eqnarray}
and
\begin{eqnarray}\label{exam203-8}
\overline{F_0}(x)&=&\textbf{1}(x<
a_0)+\sum\limits_{n=0}^{\infty}\Big(\big(\sum\limits_{i=n}^{\infty}
Ca_i^{-\alpha}-Ca_n^{-\alpha-1}(x-a_n)\big)
\textbf{1}\big(x\in[a_n,2a_n)\big)\nonumber\\
& &
+\sum\limits_{i=n+1}^{\infty}Ca_i^{-\alpha}\textbf{1}\big(x\in[2a_n,
a_{n+1})\big)\Big).
\end{eqnarray}

According to (\ref{exam203-8}), we have
\begin{eqnarray*}
\int_{0}^{\infty}\overline {F_0}(y)dy\leq a_1+\sum_{i=1}^{\infty}\int_{a_i}^{a_{i+1}}\overline {F_0}(y)dy\leq
a_1+\sum_{i=1}^{\infty}(a_i\overline {F_0}(a_i)+a_{i+1}\overline {F_0}(a_{i+1}))<\infty.
\end{eqnarray*}
Hence $\mu(F_0)<\infty$, which implies $m(F)<\infty$.
Then by $\overline {F_0}(2a_n-1)/\overline {F_0}(2a_n)=2>1$, $n\to\infty$, we know $F_0\notin\mathcal{L}$. Next we prove $F_0\in\mathcal{OS}$, it suffices to prove
$$W(x)=\int_{[x/2,x]}\overline {F_0}(x-y)F_0(dy)=O({\overline F_0}(x)).$$
When $x\in[a_n,4a_n)$ for all $n = 0,1,\cdot \cdot \cdot$, $$W(x)/\overline F_0(x)\leq\max_{y\in(a_n/2,4a_n)}f(y)\int_{[x/2,x]}\overline F_0(x-y)dy/\overline F_0(4a_n)\leq\mu(F_0)<\infty.$$
When $x\in[4a_n,a_{n+1})$ for all $n = 0,1,\cdot \cdot \cdot$, $W(x)/\overline F_0(x)=0<\infty$. Therefore, we can get $F_0^{*k}\in\cal{OS}$ for all $k\ge1$. Further, by Corollary \ref{101}, $F_0^{*k}\notin\mathcal{L}$ for all $k\ge1$.

Clearly, $F\notin\cal{L}(\gamma)$ and $F$ satisfies condition (\ref{104}). Further, by Proposition \ref{pron101} (1), we know $F\not\in\cal{OS}$. However, condition (\ref{1040}) does not hold for this $F$. In fact, by (\ref{3036}), $F_0\in\mathcal{OS}$ and $\mu(F_0)<\infty$,
\begin{eqnarray*}
\liminf_{n\to\infty}\overline F(4a_n)/\overline {F^{*2}}(4a_n)&=&\liminf_{n\to\infty}\overline F_0(4a_n)/\Big(\overline {F_0^{*2}}(4a_n)+2\gamma\int_{[0,2a_n]}\overline {F_0}(4a_n-y)\overline {F_0}(y)dy\Big)\nonumber\\
&\ge& 1/(C^*(F_0)+2\gamma\mu(F_0))>0.
\end{eqnarray*}
For the conclusion that $F^{*k},k\ge2$ and $F^{*\tau}$ belongs to the class $\mathcal{L}(\gamma)$,
by Lemma \ref{lemma202}, Corollary \ref{Corol201} (1), Theorem \ref{thm203} and Corollary \ref{Corol102},
we only need to show $F^{*2}\in\mathcal{L}(\gamma)$.
To this end, similar to the proof of Proposition \ref{pron301}, we just need to prove (\ref{exam301-70}).

For all $n = 0,1,\cdot \cdot \cdot$, because $W(x)=0$ for $x\in[4a_n,a_{n+1})$, we only need to deal with $W(x)$ in the following two cases:
$i)\ x\in[a_n,3a_n)$ and $ii)\ x\in[3a_n,4a_n)$.

In case $i)$, by (\ref{exam203-7}), (\ref{exam203-8}) and $\overline{F_0}((x+2^{-1}a_n)/2)\ge\overline{F_0}(7a_n/4)$, we have,
\begin{eqnarray}\label{exam203-9}
W(x)/T(x)&=&\int_{[a_n,x]}\overline{F_0}(x-y)F_0(dy)/
\int_{[x/2,(x+2^{-1}a_n)/2]}\overline{F_0}(y)\overline{F_0}(x-y)dy\nonumber\\
&\leq&4Ca_n^{-\alpha-2}\int_{[0,x-a_n]}\overline{F_0}(y)dy/
\lo(\overline{F_0}(x/2)\overline{F_0}((x+2^{-1}a_n)/2)\ro)\nonumber\\
&\leq&4\mu(F_0)Ca_n^{-\alpha-2}/
\lo(\overline{F_0}(3a_n/2)\overline{F_0}(7a_n/4)\ro)\nonumber\\
&=&O(a_n^{\alpha-2})\to0,~~~n\rightarrow\infty.
\end{eqnarray}

Then, in case $ii)$, by (\ref{exam203-7}), (\ref{exam203-8}) and $\int_{x/2}^{x}\overline{F_0}(x-y)\overline{F_0}(y)dy\ge
Ca_n^{-\alpha-1}\int_{0}^{x/2}\overline{F_0}(y)dy$. Therefore
\begin{eqnarray}\label{exam203-10}
W(x)/T(x)&\leq&
\int_{[x/2,2a_n]}\overline{F_0}(x-y)F_0(dy)/
\Big(\overline{F_0}(a_{n+1})\int_{[0,x/2]}\overline{F_0}(y)dy\Big)\nonumber\\
&\leq&\int_{[3a_n/2,2a_n]}\overline{F_0}(3a_n-y)dy/
\Big(\int_{[0,3a_n/2]}\overline{F_0}(y)dy\Big)\nonumber\\
&=&O(a_n^{1-\alpha})\to0,~~~n\rightarrow\infty.
\end{eqnarray}

According to (\ref{exam203-9}) and (\ref{exam203-10}), (\ref{exam301-70}) holds,
thus $F^{*2}\in\cal{L}(\gamma)$ by Lemma \ref{lemma302}.

In the following, we go to prove $F^{*k}\in\mathcal{OS}, \ k \ge 2.$ By Proposition 2.6 of \cite{SW2005-1} we only need to prove $F^{*2}\in\mathcal{OS},$
%By $F^{*2}\in\cal{L}(\gamma)$ and $m(F)<\infty$ it suffices to prove
{or equivalently, to prove that
\begin{eqnarray}\label{40311}
%\lim_{A\to\infty}\limsup
\int_{[x/2,x]}\overline {F^{*2}}(x-y)F^{*2}(dy)=O\big(\overline {F^{*2}}(x)\big).
\end{eqnarray}
%From $W(x)=O(T(x))$, it is easy to verify $\overline {F^{*2}}(x)=O(e^{-\gamma x}T(x))$. Since $F_0$ is an absolutely continuous distribution, then $F_0^{*2}$ and $F^{*2}$ have density function,
We denote the density of $F_0^{*2}$ and $F^{*2}$ by $f_0^{\otimes2}$ and $f^{\otimes2}$ respectively.} When $x\in[a_n,a_{n+1})$, by (\ref{exam203-7}) and (\ref{exam203-8}) , we have
\begin{eqnarray*}
f_0^{\otimes2}(x)=2\int_{[x/2,x]}f_0(x-y)f_0(y)dy\leq 2Ca_n^{-\alpha-1}=O(\overline {F_0}(x)).
\end{eqnarray*}
Thus by (\ref{3036}) and $W(x)=O(T(x))$ we know
\begin{eqnarray*}
f^{\otimes2}(x)=\gamma\overline {F^{*2}}(x)+e^{-\gamma x}\lo(f_0^{\otimes2}(x)+\int_{[0,x]}\overline {F_0}(x-y)F_0(dy)-\overline {F_0}(x)\ro)={O(e^{-\gamma x}T(x)).}
\end{eqnarray*}
Hence, {in order to prove (\ref{40311}),} we only need to prove
\begin{eqnarray}\label{40312}
R(x)=\int_{[x/2,x]}T(x-y)T(y)dy=O(T(x)).
\end{eqnarray}
By (\ref{exam203-8}),
\begin{eqnarray}\label{6001}
T(x)\leq \overline {F_0}(a_i)\mu(F_0)\ \ {for\ all\ x\in[4a_{i-1},4a_i),\ i\ge1.}
\end{eqnarray}
{Further, we know that}
\begin{eqnarray}\label{6002}
T(x)\leq 16C\mu(F_0){x}^{-\alpha}\ \ {for\ all}\ x>4a_1,
\end{eqnarray}
thus
\begin{eqnarray*}
\mu(T)=\int_{[0,\infty)}T(y)dy
{<\infty}.
\end{eqnarray*}
Furthermore, for all $a_0<x_1<x_2<\infty$, $x\in[x_1,x_2]$, we have
\begin{eqnarray}\label{60001}
T(x)&\leq&\lo(\int_{[x_1/2,x_1]}+\int_{(x_1,x]}\ro)
\overline{F_0}(x-y)\overline{F_0}(y)dy\nonumber\\
&\leq& T(x_1)+\overline{F_0}(x_1)\mu(F_0)\nonumber\\
&\leq& (1+\mu(F_0)a_0^{-1})T(x_1).
\end{eqnarray}

{Now}, we deal with $R(x)$ in the following four cases:
$i)\ x\in[a_n,3a_n/2)$; $ii)\ x\in[3a_n/2,3a_n)$; $iii)\ x\in[3a_n,4a_n-4a_n^{2-\alpha})$;
and $iv)\ x\in[4a_n-4a_n^{2-\alpha},a_{n+1})$.

In case $i)$, by (\ref{exam203-8}), (\ref{6001}) and
$T(x)\geq \overline{F_0}(3a_n/2)\int_{[0,x/2]}\overline{F_0}(y)dy,$
we have
\begin{eqnarray}\label{6003}
R(x)/T(x)&\leq&\overline {F_0}(a_n)\mu(F_0)\int_{[0,3a_n/4]}
{T}(y)dy
/\overline {F_0}(3a_n/2)\int_{[0,a_n/2]}
\overline{F_0}(y)dy\nonumber\\
&
{\to}&
{2\mu(T)}<\infty,~~~n\rightarrow\infty.
\end{eqnarray}

Then, in case $ii)$, by (\ref{exam203-8}), (\ref{6001}), (\ref{6002}), (\ref{60001}),
$$T(x)\geq \int_{[3a_n/2,7a_n/4]}\overline {F_0}(3a_n-y)\overline {F_0}(y)dy \geq C^2a_n^{1-2\alpha}/2^{5},$$ $$T(x-4a_{n-1})/T(x)=O(a_{n-1}a_n^{\alpha-2})\ \ and\ \
T(x-a_n^{2-\alpha})=O(T(x)).$$

Thus,
\begin{eqnarray}\label{6004}
&&R(x)/T(x)=\Big(\int_{[x/2,x-4a_{n-1}]}+\int_{(x-4a_{n-1},x-a_n^{2-\alpha}]}+\int_{(x-a_n^{2-\alpha},x]}\Big)T(x-y)T(y)dy/T(x)\nonumber\\
&=&O\Big(\int_{[x/2,x-4a_{n-1}]}T(x-y)T(y)dy
+T(x-4a_{n-1})\int_{[a_n^{2-\alpha},4a_{n-1}]}T(y)dy\nonumber\\
& &\ \ \ \ \ \ \ \ \ +T(x-a_n^{2-\alpha})\int_{[0,a_n^{2-\alpha}]}T(y)dy\Big)/T(x)\nonumber\\
&=&{O\Big(\overline {F_0}(a_n)\int_{[4a_{n-1},3a_n/2]}T(y)dy/a_n^{1-2\alpha}
+a_{n-1}a_n^{\alpha-2}\int_{[a_n^{2-\alpha},4a_{n-1}]}y^{-\alpha}dy+\int_{[0,a_n^{2-\alpha}]}}T(y)dy\Big)\nonumber\\
&=&O\big(a_n^{(\alpha+1)^{-1}(\alpha^3-\alpha^2-\alpha)}+1\big)<\infty,~~~n\rightarrow\infty.
\end{eqnarray}

Next, in case $iii)$, by (\ref{exam203-8})
and (\ref{6001}),
we know that, when $y\in[x/2,2a_n-2a_n^{2-\alpha})$, $x-y\in[x-2a_n+2a_n^{2-\alpha},x/2)$,
$$T(y)=O(\overline{F_0}(y)),\ \ \ \ T(x-y)=O(\overline{F_0}(x-y));$$
and
$$T(x-4a_{n-1})=O(T(x)),\ \ \ \ T(x)\ge\overline {F_0}(2a_n)\int_{[0,x-2a_n]}\overline{F_0}(y)dy.$$
Thus, by (\ref{60001}),
\begin{eqnarray}\label{6005}
R(x)/T(x)&=&\Big(\int_{[x/2,2a_n-2a_n^{2-\alpha}]}+\int_{(2a_n-2a_n^{2-\alpha},x-4a_{n-1}]}+\int_{(x-4a_{n-1},x]}\Big)T(x-y)T(y)dy
/T(x)\nonumber\\
&=&O\Big(\int_{[x/2,2a_n-2a_n^{2-\alpha}]}\overline{F_0}(x-y)\overline{F_0}(y)dy/T(x)
+\int_{[0,4a_{n-1}]}T(y)dy\nonumber\\
& &\ \ \ \ \ \ \ \ \ +a_n T(2a_n-2a_n^{2-\alpha})T(4a_{n-1})/\overline {F_0}(2a_n)\int_{[0,a_n]}\overline{F_0}(y)dy\Big)\nonumber\\
&=&O(a_n^{3-2\alpha}+1)<\infty,~~~n\rightarrow\infty.
\end{eqnarray}

For case $iv)$,by (\ref{exam203-8}), (\ref{6001}), (\ref{60001}),
and when $x\in[4a_n-4a_n^{2-\alpha},a_{n+1})$,
$$T(x)\geq\overline {F_0}(2a_n)\int_{[0,x-2a_n]}\overline{F_0}(y)dy.$$
Thus
\begin{eqnarray}\label{6006}
&&R(x)/T(x)\leq\Big(\int_{[2a_n-2a_n^{2-\alpha},4a_n-8a_{n-1})}+\int_{[4a_n-8a_{n-1},x]}\Big)T(x-y)T(y)dy
\int_{[0,a_n]}\overline{F_0}(y)dy/\overline {F_0}(2a_n)\nonumber\\
&=&O(a_n T(2a_n-2a_n^{2-\alpha})T(4a_{n-1})+T(4a_n-8a_{n-1})\mu(T))
\int_{[0,a_n]}\overline{F_0}(y)dy/\overline {F_0}(2a_n)\nonumber\\
&=&O(a_n^{2-3\alpha}+a_n^{-\alpha-1})
\int_{[0,a_n]}\overline{F_0}(y)dy/\overline {F_0}(2a_n)\nonumber\\
&=&O\big(a_n^{3-2\alpha}+1\big)<\infty,~~~n\rightarrow\infty.
\end{eqnarray}
By (\ref{6003})-(\ref{6006}), we can get (\ref{40312}) holds. Thus $F^{*2}\in\mathcal{OS}$.

{Finally, from the following formula, we know $F^{*2}\notin\mathcal{S}(\gamma)$:
\begin{eqnarray*}
&&\int_{[a_n,2a_n]}\overline {F^{*2}}(2a_n-y)F^{*2}(dy)\ge\int_{[a_n,3a_n/2]}\overline {F^{*2}}(2a_n-y)f_0^{\otimes2}(y)dy\nonumber\\
&\ge&e^{-2\gamma a_n }\int_{[a_n,3a_n/2]}T(2a_n-y)T(y)dy\approx\overline {F^{*2}}(2a_n)
\end{eqnarray*}
where $g_1(x)\approx g_2(x)$ means that $g_1(x)=O(g_2(x))$ and $g_2(x)=O(g_1(x))$ for two positive functions $g_1$ and $g_2$. Then by $F^{*k}\in\mathcal{L}(\gamma)\cap\mathcal{OS}$ for all $k\ge2$, we know $\overline {F^{*2}}(x)\approx \overline {F^{*k}}(x)$, which is called that distribution $F^{*k}$ is weakly tail equivalent to distribution $F^{*2}$ for all $k\ge2$. By Lemma 2.6 in \cite{W2008} and $F\notin\mathcal{S}(\gamma)$, we can get $F^{*k}\notin\mathcal{S}(\gamma)$ for all $k\ge3$.}
\hfill$\Box$
\vspace{0.2cm}

Finally, we introduce a distribution which fails to satisfy the condition (\ref{104}), thus it is different from the distributions in the classes $\mathcal{F}_i(\gamma),i=1,2,3$.
\begin{defin}\label{defin302}
{For any $\gamma>0$, denote}
$$\mathcal{F}_4(\gamma)=\{F,\ \overline F(x)=\textbf{\emph{1}}(x< 0)+\sum_{k=0}^{\infty}(e^{-\gamma x}+e^{-\gamma e^{k+1}})\textbf{\emph{1}}(e^{k}\le x<e^{k+1})/2\ \ for\ all\ x\}.$$
\end{defin}

The distribution family $\mathcal{F}_4(\gamma)$ was introduced by Theorem 3.1 {in} \cite{SW2005-2} with properties {that $F\notin \mathcal{L}(\gamma)$ and $F^{*2}\in \mathcal{L}(\gamma)$. Here, we give some of its new properties.}
\begin{pron}\label{pron304}
If $F\in\mathcal{F}_4(\gamma)$ for some $\gamma>0$, then $F^{*k}\in \mathcal{L}(\gamma )$ for all $k\ge 2$. {In addition, for the distribution $F$,  condition (\ref{1040}) is satisfied, but condition (\ref{104}) is not}. Further, if the condition (\ref{21}) is satisfied with some non-negative integer-valued random variable $\tau$, then $F^{*\tau}\in \mathcal{L}(\gamma)$, while $F\in\cal{OL}\setminus\mathcal{L}(\gamma)$.
\end{pron}

\proof From Theorem 3.1 in \cite{SW2005-2}, we know that $F\notin\mathcal{L}(\gamma)$ and $F^{*2}\in \mathcal{L}(\gamma )$. Further, $F\in\cal{OL}$ is an obvious fact, and $F^{*k}\in \mathcal{L}(\gamma )$ for all $k\ge 3$ follows from Remark \ref{rem203} of the present paper.
{In addition, it is easy to verify that $\overline {F^{*2}}(x)\sim \overline {\mu^{*2}}(x)/4$, where $\mu$ is a standard exponential distribution such that $\overline {\mu}(x)=\textbf{1}(x<0)+e^{-\gamma x}\textbf{1}(x\ge0)$ for all $x$. Hence, condition (\ref{1040}) is satisfied. However, condition (\ref{104}) does not hold, in fact, for any $t>0$,
$$\liminf_{k\to\infty}\overline F(e^{k+1}-2t)/\overline F(e^{k+1}-t)=(e^{2\gamma t}+1)/(e^{\gamma t}+1)<e^{\gamma t}.$$}

{Finally, $F^{*\tau}\in \mathcal{L}(\gamma)$ follows from (\ref{1040}), (\ref{21}) and Theorem \ref{thm203}.}\hfill$\Box$
\vspace{0.2cm}

\begin{remark}\label{rem401}
Here we give a distinction between distributions in $\mathcal{F}_i(\gamma),i=1,2,3,4$ and in Lemma 3.1 of \cite{W2015}. For some $\gamma>0$ and $i=1,2,3$, it is easy to verify that, if $F\in\mathcal{F}_i(\gamma)$, then
$$\limsup_{t\to\infty}\limsup\overline{F}(x-t)\big(\overline{F}(x)\big)^{-1}=\infty.$$
Therefore, the distribution $F$ is not weakly tail equivalent to any distribution in $\mathcal{L}(\gamma)$. However, the distribution in $\mathcal{F}_4(\gamma)$ is weakly tail equivalent to a standard exponential distribution, and the distribution in \cite{W2015} is weakly tail equivalent to a distribution in the class $\mathcal{S}(\gamma)$.
\end{remark}

Now, we prove Theorem \ref{thm101}. For $i=1,2,3,4$, we take any distribution $F\in\mathcal{F}_i(\gamma)$ and distribution $H_2$ such as in (\ref{101}), then by Propositions \ref{pron301}, \ref{pron302}, \ref{pron303} or \ref{pron304}, we have $H_2\in\mathcal{L}(\gamma)$. Further, we take distribution $H_1$ satisfying (\ref{103}), or more simply, we take a distribution $H_1$ such that $\overline{H_1}=O(e^{-\beta x})$ for any $\gamma>\beta>0$. Thus, by Theorem \ref{thm201} of the paper or Theorem 3 in \cite{EG1980}, the infinitely divisible distribution $H=H_1*H_2\in\mathcal{L}(\gamma)$.

{In particular, for some $\gamma>0$, we choose any distribution $F\in\mathcal{F}_3(\gamma)$ and distribution $H_1$ such that $\overline{H_1}(x)=o(\overline{H_2}(x))$, then by Proposition \ref{pron303}, Theorem \ref{thm203} and Theorem \ref{thm204}, the infinitely divisible distribution $H\in \big(\mathcal{L}(\gamma )\cap\mathcal{OS}\big)\setminus\mathcal{S}(\gamma)$ and the inequality (\ref{10400}) holds.}

So far, we have completed all the proof of Theorem \ref{thm101}.

\section{Some remarks}
\setcounter{thm}{0}\setcounter{Corol}{0}\setcounter{lemma}{0}\setcounter{pron}{0}\setcounter{equation}{0}
\setcounter{remark}{0}\setcounter{exam}{0}\setcounter{property}{0}\setcounter{defin}{0}

In this Section, we first note that condition (\ref{103}) can not be deduced by condition (\ref{102}). Then, we provide examples showing that conditions (\ref{l30101}) and (\ref{l3011}) cannot be deduced from each other. Finally, we give a local version on the study of Embrechts-Goldie's
conjecture.

\subsection{On (\ref{102}) and (\ref{103})}

In Example 5.1 of \cite{XFW2015}, there are two distributions $F_{1}$ and $F_{2}$ such that $F_{1}\in\mathcal{OL}\backslash\mathcal{L}$, $F_{2}\in\mathcal{S}$, $\overline{F_{2}}(x)=o(\overline{F_{1}}(x))$, $F_{1}*F_{2}\in\mathcal{L}$ and
$\overline{F_{1}}(x-t)-\overline{F_{1}}(x)=o(\overline{F_{2}}(x))$.

Then we take $F_{i0}=F_i,i=1,2$. As (\ref{303}), for any constant $\gamma>0$ and distribution $F_{i0},i=1,2$, we define the distribution $F_{i\gamma}$ such that
\begin{eqnarray*}
\overline{F_{i\gamma}}(x)=\textbf{\emph{\emph{1}}}(x<0)+e^{-\gamma x}\overline{F_{i0}}(x)\textbf{\emph{\emph{1}}}(x\ge0),\ x\in(-\infty,\infty).
\end{eqnarray*}
Then $F_{1\gamma}\in\mathcal{OL}\backslash\mathcal{L}(\gamma)$, $F_{2\gamma}\in\mathcal{L}(\gamma)$, $\overline{F_{2\gamma}}(x)=o(\overline{F_{1\gamma}}(x))$, $\overline{F_{1\gamma}}(x-t)-e^{\gamma t}\overline{F_{1\gamma}}(x)=o(\overline{F_{2\gamma}}(x))$ and $F_{1\gamma}*F_{2\gamma}\in\mathcal{L}(\gamma)$ which comes from Theorem \ref{thm201} and the properties mentioned above.

Now, we take $H_1=F_{1\gamma}$ and
\begin{eqnarray*}
H_2(x)=e^{-\mu}\sum_{k=0}^\infty F_{2\gamma}^{*k}(x)\mu^k/k!,\ x\in(-\infty,\infty).
\end{eqnarray*}
Then, there is a positive constant $a$ such that $\overline{H_2}(x)\sim a\overline{F_{2\gamma}}(x)$. Thus, (\ref{103}) holds, while (\ref{102}) doesn't.

\subsection{Comparison of conditions (\ref{l30101}) and (\ref{l3011})}
Firstly, let $X$ be a random variable with a distribution $F_0$ in Example 3.3 of $\cite{XSWC2015}$ such that
\begin{eqnarray*}
&&\overline{F_0}(x)=\textbf{1}(x<0)+\big(x_1^{-1}(x_1^{-\alpha}-1)x+1\big)\textbf{1}(0\leq x< x_1)\nonumber\\
&+&\sum\limits_{n=1}^{\infty}\Big(\big(x_n^{-\alpha}+(x_n^{-\alpha-2}
-x_n^{-\alpha-1})(x-x_n)\big)\textbf{1}(x_n\leq x<2x_n)
+x_n^{-\alpha-1}\textbf{1}(2x_n\leq x< x_{n+1})\Big)
\end{eqnarray*}
for all $x\in(-\infty,\infty)$, where $\alpha\in(5,\infty)$, $x_1>4^{\alpha}$ and for all
integers $n\geq1$, $x_{n+1}=x_n^{1+{\alpha}^{-1}}, n = 0,1,\cdot \cdot \cdot$.
We already know $F\notin\mathcal{L}$ and $\mu_2=EX^2<\infty$. Further, we denote the mean of $X$ by $\mu_1$ and the density of $F_0$ by $f_0$. For all $s>0$ and $n = 0,1,\cdot \cdot \cdot$, when $x\in[x_n,2x_n+s)$,
$$\overline {F_0}(x-s)-\overline {F_0}(x)\leq sx_n^{-\alpha-1}=O\big(f_0(x-s)+f_0(x)\big);$$
when $x\in[2x_n+s,x_{n+1})$,
$$\overline {F_0}(x-s)-\overline {F_0}(x)=f_0(x-s)+f_0(x)=0.$$
Hence, condition (\ref{l30101}) holds for the distribution $F_0$.
Also denote $W(x)=\int_{[x/2,x]}\overline {F_0}(x-y)F_0(dy)$ and $T(x)=\int_{[x/2,x]}\overline{F_0}(x-y)\overline{F_0}(y)dy$, however,

\begin{eqnarray*}
W(2x_n)/T(2x_n)&\geq&\int_{[0,x_n]}\overline {F_0}(y)dy/\int_{[0,x_n]}\overline{F_0}(y)(1+y)dy\nonumber\\
&\to&\mu_1/(\mu_1+2^{-1}\mu_2)>0,~~~n\to\infty,
\end{eqnarray*}
which implies that $F_0$ does not satisfy condition (\ref{l3011}).

Then, let the distribution $F_1$ follow
\begin{eqnarray*}
\overline {F_1}(x)=\textbf{1}(x< 1)+\sum_{n=1}^{\infty}\Big(\textbf{1}(2n-1\le x < 2n)/(2x-2n+1)+\textbf{1}(2n\le x < 2n+1)/(2n+1)\Big)
\end{eqnarray*}
with its density
$$f_1(x)=2\sum_{n=1}^{\infty}\textbf{1}(2n-1\le x < 2n)/(2x-2n+1)^2,$$
for all $x\in(-\infty,\infty)$. By
$$\overline{F_1}(2n-3\cdot2^{-1})-\overline{F_1}(2n+2^{-1})=3/((2n-2)(2n+1))$$
and
$$f_1(2n-3\cdot2^{-1})=f_1(2n+2^{-1})=0,$$
we know that $F$ does not satisfy condition (\ref{l30101}).
On the other hand, it is also easy to verify that $f_1(x)=o(\overline {F_1}(x)),$
hence $F_1\in\mathcal{L}$ and satisfies condition (\ref{l3011}).

However, we need such a distribution which does not belong to $\mathcal{L}$. To this end, let $y_0\ge0$ and $a>1$ be two constants such that $a\overline{F_1}(y_0)\le1$. As Example 3.1 {in Xu et al.} \cite{XSWC2015}, we define a new distribution $F_0$ by
\begin{eqnarray*}
\overline{F_0}(x)&=&\overline{F_1}(x)\textbf{1}(x<x_1)+
\sum_{i=1}^{\infty}\Big(\overline{F_1}(x_i)\textbf{1}(x_i\le x<y_i)+\overline{F_1}(x)\textbf{1}(y_i\le x<x_{i+1})\Big)
\end{eqnarray*}
for $x\in (-\infty,\infty)$, where $\{x_i,i\ge1\}$ and $\{y_i, i\ge 1\}$ are two sequences of positive constants satisfying $x_i<y_i<x_{i+1}$, $\overline{F_1}(x_i)=a\overline{F}_1(y_i),\ i\ge1,\ y_i-x_i\to\infty$ and $x_{i+1}-y_i\to\infty$ as $i\to\infty$. It is easy to see that  $\overline{F_1}(x)\le\overline {F_0}(x)\le a\overline{F_1}(x)$ and $\lim_{n \to \infty}\overline {F_0}(y_n-1)/\overline {F_0}(y_n)=a>1$. Thus, the distribution $F_0$ satisfies condition (\ref{l3011}) and does not belong to class $\mathcal{L}$, while similar to $F_1$, $F_0$ does not satisfy condition (\ref{l30101}).

Therefore, the conditions (\ref{l30101}) and (\ref{l3011}) can not be deduced from each other.

\subsection{A local version of Embrechts-Goldie's conjecture}
We say that a distribution $F$ belongs
to the distribution class $\mathcal{L}_{loc}$,
if for all $T>0$ and $t \neq 0$, $F(x+\Delta_T)>0$ for all $x$ large enough and
$$F(x-t+\Delta_T)=F(x-t,x-t+T]\sim F(x,x+T]=F(x+\Delta_T).$$
Further, if a distribution $F$ belongs to the class $\mathcal{L}_{loc}$,
and if for all $t \neq 0$ and for all $T>0$,
$$F^{*2}(x+\Delta_T)\sim2 F(x+\Delta_T),$$
then we say that the distribution $F$ belongs to the distribution class $\mathcal{S}_{loc}$. The concepts of the classes $\mathcal{L}_{loc}$ and $\mathcal{S}_{loc}$ can be found in Borovkov and Borovkov \cite{BB2008}. Similar to the classes $\mathcal{OL}$ and $\mathcal{OS}$, we can also give the concepts of the classes $\mathcal{OL}_{loc}$ and $\mathcal{OS}_{loc}$, respectively.

According to Proposition 2.1 of Wang and Wang \cite{WW2011} and Proposition \ref{pron303} of the paper, the following conclusion holds.
\begin{pron}\label{pron601}
The distribution class $\big(\mathcal{L}_{loc}\cap\mathcal{OS}_{loc}\big)\setminus\mathcal{S}_{loc}$ is not closed under convolution roots.
\end{pron}

Corollary 1.1 (i) in \cite{W2015} notes that the distribution class $\mathcal{S}_{loc}$ is not closed under convolution roots too. The conclusion with the Proposition \ref{pron601} of the paper constitutes a more complete answer to the Embrechts-Goldie's conjecture in local sense. A further discussion on the topic for the local distribution class can be found in Wang et al. \cite{WXC2016}.
\\
\\
\textbf{Acknowledgements} The authors are very grateful to Sergey Foss and Vsevolod Shneer for their helpful discussions and suggestions. The authors are also very grateful to Jonas $\check{S}$iaulys for his detailed and valuable comments and suggestions for the whole content of this paper.

\vspace{0.2cm}

\end{document}